\documentclass[12pt,final]{article}
\usepackage{amsmath,ascmac,latexsym,fancybox,showkeys,color}
\usepackage[pagewise]{lineno}

\def\R{{\bf R}}

\def\de{{\partial}}
\def\dy{\displaystyle}
\def\ep{{\varepsilon}}
\def\w{\widetilde}
\def\-{|\!|\!|}
\def\supp{\mbox{supp }}

\newtheorem{thm}{Theorem}
\newtheorem{lm}{Lemma}[section]
\newtheorem{prop}{Proposition}[section]

\newtheorem{re}{Remark}[section]
\makeatletter
 	
 \@addtoreset{equation}{section}
\makeatother

\title{The lifespan estimates of radially symmetric solutions to 
systems of nonlinear wave equations in even space dimensions}
\author{Yuki Kurokawa
\thanks{Hokkaido University of Education, Kushiro 085-8580, Japan. 
e-mail : kurokawa.yuki@k.hokkyodai.ac.jp.}}

\date{{\small $
\begin{array}{l}
\mbox{{\it MSC 2020 :} 35L71, 35L05, 35E15}\\
\mbox{{\it Key words :} system, semilinear wave equation, radially symmetric solution, lifespan}\\
\mbox{{\it Running head :} Lifespan estimates of radially symmetric solutions}
\end{array}$
}}

\begin{document}
\maketitle
\begin{abstract}
\par
The optimal lifespan estimates of a solution to weakly coupled systems of wave equations have been investigated by many works, except for the lower bound in even space dimensions. 
Our aim is to prove the open part under the assumption of radial symmetry on the solution. 
The odd dimensional case was already obtained in our previous paper \cite{Kodd} by long time existence of the solution in weighted $L^\infty$ space. 
In this paper, we employ similar methods. 
The main difficulty is found in estimating the integral kernel which is completely different from odd dimensional case. 
\end{abstract}
	\section{Introduction}
\label{sec:intro}
\par
We consider the following Cauchy problem of the systems 
of semilinear wave equations :
\begin{equation}
\left\{
\begin{array}{l}
u_{tt}-\Delta u=|v|^p,\\
v_{tt}-\Delta v=|u|^q,
\end{array}
\right.
\qquad\mbox{in\quad}\R^n\times (0,\infty),
\label{D}
\end{equation}
for $p,q>1$ with the initial data
\begin{equation}
\left\{
\begin{array}{l}
u(x,0)=\ep \varphi_1(x),\ u_t(x,0)=\ep \psi_1(x),\\
v(x,0)=\ep \varphi_2(x),\ v_t(x,0)=\ep \psi_2(x),
\end{array}
\right.\qquad\mbox{for}\quad x\in\R^n, 
\label{ID}
\end{equation}
where $\ep >0$ is a sufficiently small parameter and $\varphi_i$, $\psi_i$ ($i=1,2$), which present the shape of the data, 
are sufficiently smooth and have compact support. 
\par 
For this problem, Del Santo, Georgiev and Mitidieri \cite{DGM} proved that 
\begin{itemize}
\item $F(p,q;n)<0$, $2\le p,q\le 3$ $\Longrightarrow$ Existence of a global-in-time solution for $\ep<<1$,
\item $F(p,q;n)>0$, some positivety on data $\Longrightarrow$ Nonexistence of of a global-in-time solution, 
\end{itemize}
for $n\ge 2$, where 
\begin{equation}
F(p,q;n):=\max\left\{
\dy\frac{p+2+q^{-1}}{pq-1}, \dy\frac{q+2+p^{-1}}{pq-1}
\right\}
-\dy\frac{n-1}{2}.
\label{F}
\end{equation}
Their restrictions ``$p,q\le 3$'' and ``some positivity on data'' were relaxed by Del Santo \cite{Del} in three space dimensions. 
Thus, we see that there exists a critical curve $F(p,q;n)=0$ on the $p$-$q$ plane, which is symmetric across the line $p=q$. 
Independently, Deng \cite{D97}, \cite{D99} proved the blow-up of a weak solution for the case of $F(p,q;n)>0$ with $p\le q$ and $n\ge 2$ for some positive data. 
The blow-up of a solution for the critical case $F(p,q;n)=0$ was proved by Del Santo and Mitidieri \cite{DM} in three space dimensions under some positivety on data. 
\par 
As for the blow-up case {\it i.e.} $F(p,q)\ge 0$, there are several results on "lifespan", the maximal existence time of a solution. 
Kubo and Ohta \cite{KO} obtained the upper bound of the lifespan of the classical solution in two or three space dimensions for some positive data. 
The optimal lifespan of a classical solution for $n=2$ or $3$ has been proved for all cases by 
Agemi, Kurokawa and Takamura \cite{AKT} and Kurokawa and Takamura \cite{KT2}. 
Thus, we have 
\begin{align*}
&\quad\ T_0(\ep)=\infty &\mbox{if}\ &F(p,q;n)<0,\\
\exp(c\ep^{-\min\{p(pq-1),q(pq-1)\}})&\le T_0(\ep)\le\exp(c'\ep^{-\min\{p(pq-1),q(pq-1)\}}) &\mbox{if}\ &F(p,q;n)=0,\ p\neq q,\\
\exp(c\ep^{-p(p-1)\}})&\le T_0(\ep)\le\exp(c'\ep^{-p(p-1)}) &\mbox{if}\ &F(p,q;n)=0,\ p=q,\\
c\ep^{-F(p,q;n)^{-1}}&\le T_0(\ep)\le c'\ep^{-F(p,q;n)^{-1}} &\mbox{if}\ &F(p,q;n)>0,
\end{align*}
for $n=2,3$, where $c$ and $c'$ are positive constants independent of $\ep$ and $T_0(\ep)$ means the lifespan of the classical solution to (\ref{D}) with (\ref{ID}). 

\par
The higher dimensional case was first investigated by Georgiev, Takamura and Zhou in \cite{GTZ} for the sub-critical case. 
They showed that for some data of Sobolev functions the lifespan $T_{\ep}$ of a certain weak solution in Lebesgue space  can be pintched from both side as follows. 
\begin{equation}
\label{error1}
C_1\ep^{-(F(p,q;n)+\sigma)^{-1}}\le T_{\ep}\le C_2\ep^{-F(p,q;n)^{-1}},
\end{equation}
where $C_1$, $C_2$ are positive constants independent of $\ep$ and 
$\sigma>0$ is a small error. 
Moreover, Remark 1.2 in \cite{GTZ}, they said that they can get the similar lower bound of the lifespan for the case of $F(p,q;n)=0$ :
\begin{equation}
\label{error2}
C_1\ep^{-1/\sigma}\le T_\ep. 
\end{equation}
Our previous paper, Kurokawa \cite{Kodd}, succeeded to remove the error $\sigma$ in the lower bound of the lifespan for $F(p,q;n)\ge 0$ in odd space dimensions under the assumption of radial symmetricity of a solution. 
The present paper is devoted to remove the error for a radially symmetric solution in even space dimensions as well. 
\par
The upper bound of the lifespan for the critical case $F(p,q;n)=0$ was obtained by Kurokawa, Takamura and Wakasa \cite{KTW} for $n\ge 4$. 
Their estimate is same as the one in low space dimensions. 
Recently there are some alternative proof for the upper bound of the lifespan. 
Ikeda, Sobajima and Wakasa \cite{ISW} introduced new test functions in the weak form to obtain the same results for $n\ge 2$ and $F\ge 0$. 
Palmieri and Takamura \cite{PT} studied the systems of semilinear damped wave equations : 
\[
\left\{
\begin{array}{ll}
u_{tt}-\Delta u+b_1(t)u_t=|v|^p,&x\in\R^n,\ t>0,\\
v_{tt}-\Delta v+b_2(t)v_t=|u|^q,&x\in\R^n,\ t>0,\\
(u,u_t,v,v_t)(0,x)=(\ep u_0,\ep u_1,\ep v_0,\ep v_1)(x)&x\in\R^n
\end{array}
\right.,
\]
for $b_1,b_2\in C([0,\infty))\cap L^1([0,\infty))$ and they get the upper bound of the lifespan of an energy solution for $n\ge 2$ and $F\ge 0$. 
This work includes the case of $b_1\equiv 0$ and $b_2\equiv 0$ {\it i.e.} our systems (\ref{D}). 
\par
From now on, we consider the following radially symmetric version of the systems of wave equations for $r=|x|$ and $T>0$ : 
\begin{equation}
\left\{
\begin{array}{l}
u_{tt}-u_{rr}-\dy\frac{n-1}{r}u_r=|v|^p,\\
v_{tt}-v_{rr}-\dy\frac{n-1}{r}v_r=|u|^q,
\end{array}
\right.
\qquad\mbox{in\quad}(0,\infty)\times (0,T)
\label{D_rad}
\end{equation}
for $p,q>1$ with the initial data
\begin{equation}
\left\{
\begin{array}{l}
u(r,0)=\ep f_1(r),\ u_t(r,0)=\ep g_1(r),\\
v(r,0)=\ep f_2(r),\ v_t(r,0)=\ep g_2(r),
\end{array}
\right.
\label{ID_rad}
\end{equation}
where, $\ep >0$ and $f_i$, $g_i$ ($i=1,2)$ are sufficiently smooth functions of compact support. 
More precisely, we assume that 
\begin{equation}
\label{supp_data}
\supp \{f_1,f_2,g_1,g_2\}\subset\{r\le t+k\}\quad\mbox{with}\quad k>1. 
\end{equation}
In this paper, we study the following system of the integral equations which is equivalent to (\ref{D_rad}) with (\ref{ID_rad}). 
\begin{equation}
\label{I_rad}
\left\{
\begin{array}{l}
u(r,t)=\ep u^0(r,t)+L(|v|^p)(r,t),\\
v(r,t)=\ep v^0(r,t)+L(|u|^q)(r,t), 
\end{array}
\right.\quad\mbox{for}\quad (r,t)\in(0,\infty)\times(0,T), 
\end{equation}
where $u^0$, $v^0$ are solutions to free wave equation with the same data as the ones of $u$, $v$ in (\ref{ID_rad}) without $\ep$. 
$L$ is an integral operator associated to Duhamel's principle which will be defined in the next section. 
In odd space dimensions, (\ref{D_rad}) with (\ref{ID_rad}) is equivalent to (\ref{I_rad}), since we can construct a $C^2$ solution of (\ref{I_rad}) (See \cite{KK_odd}, \cite{Kodd}). 
On the other hand, we do not expect to obtain $C^2$ solution due to the lack of the differentiability of the kernel of the integral expressions. 
Therefore, our solution is defined as a $C^1$ solution of (\ref{I_rad}). 
Namely, the lifespan $T(\ep)$ in the present paper is defined by 
\begin{align}
T(\ep)
&:=\sup\{T\in(0,\infty]\ :\ \mbox{There exists a radially symmetric solution}\nonumber\\
&\qquad \ (u,v)\in(C^1((0,\infty)\times[0,T]))^2\ \mbox{to the system (\ref{I_rad}) with any (\ref{ID_rad}).}\}. \label{lifespan}
\end{align}
\begin{re}
The solution defined above can be a weak solution of (\ref{D}) in the sense of $(u,v)\in (C^1((0,\infty)\times[0,T]))^2\cap (C^2([0,T];\mathcal{D}'(\R^n)))^2$. 
See the next section. More precisely, see Theorem 2.1 and Theorem 2.2 in \cite{KK_even}. 
\end{re}
\par
Now our main result is the following. 
\begin{thm}
\label{main}
Let $n=2m+2$ ($m=2,3,4,\cdots$), $F(p,q;n)\ge 0$ and 
\[
p,q\in\left(\frac{n+1}{n-1},\frac{n+3}{n-1}\right).
\] 
Assume that $f_i\in C^2([0,\infty))$ and $g_i\in C^1([0,\infty))$ ($i=1,2$) satisfy the support property (\ref{supp_data}). 
Then, there exists a positive constant $\ep_0=\ep_0(p,q,f_1,f_2,g_1,g_2,n,k)$ such that the lifespan $T(\ep)$, defined in (\ref{lifespan}), satisfies 
\[
T(\ep)\ge\left\{\begin{array}{ll}
\exp(C\ep^{-\min\{p(pq-1),q(pq-1)\}})&\mbox{if}\quad F(p,q;n)=0,\ p\neq q,\\
\exp(C\ep^{-p(p-1)})&\mbox{if}\quad F(p,q;n)=0,\ p=q,\\
C\ep^{-F(p,q;n)^{-1}}&\mbox{if}\quad F(p,q;n)>0
\end{array}\right.
\]
for any $\ep\in(0,\ep_0]$, where $C$ is a positive constant independent of $\ep$ and $F(p,q;n)$ is the one in (\ref{F}).
\end{thm}

\begin{re}
By the symmetricity of the systems (\ref{I_rad}), it is enough to prove the case of $p\le q$ only. 
\end{re}

\begin{re}
In this paper, $m=1$ {\it i.e.} four space dimensional case is excluded.  
Because we employ the methods in Kubo and Kubota \cite{KK_even} which study the asymptotic behavior of radially symmetric solutions to the corresponding single wave equations in even space dimensions with $n\ge 6$. 
Thus, the similar results for $n=4$ has not been obtained yet, while Li and Zhou \cite{LZ} and Lindblad and Sogge \cite{LS} already get the sharp lower bound for the single equation in four space dimensions. 
\end{re}

\begin{re}
\par
The corresponding single equation : 
\begin{equation}
\label{single}
\left\{
\begin{array}{ll}
u_{tt}-\Delta u=|u|^p &\mbox{in}\quad \R^n\times[0,\infty),\\
u(x,0)=\ep f(x),\quad u_t(x,0)=\ep g(x) & \mbox{for}\quad x\in\R^n
\end{array}\right.
\end{equation}
has been investigated by many works and we now know that 
the critical power $p_0(n)$ is the positive root of 
\[
\gamma(p,n)\equiv 2+(n+1)p-(n-1)p^2=0.
\]
Namely, if $p>p_0(n)$, there exists a global-in time solution to (\ref{single}) for sufficiently small $\ep$. 
On the other hand, a blow-up occurs, for some initial data, if $p\le p_0(n)$ 
(See \cite{J}, \cite{St}, \cite{Gla_G}, \cite{Gla_B}, \cite{Schaeffer}, \cite{Si}, \cite{R2}, \cite{GLS}, \cite{YZ}, \cite{Z_high}). 
\par
The lifespan estimates are also well-studied by many works. 
See \cite{Z1}, \cite{Lind}, \cite{Z3}, \cite{Z2}, \cite{T_improve}, \cite{IKTW} for low space dimensions 
and \cite{LZ}, \cite{LS}, \cite{TW}, \cite{ZH}, \cite{LaiZ} for high space dimensions. 
For the details, see \cite{IKTW} and its introduction. 
Lindblad and Sogge \cite{LS} has an important relation to our problems. 
They studied (\ref{single}) with radially symmetric data and showed that for $n\ge 3$ and $(n+1)/(n-1)\le p\le p_0(n)$ (\ref{single}) has a unique solution $u$ verifying 
\[
t^{(n-1)/2-(n+1)/2p}r^{(n+1)/2p}u\in L_t^\infty L_r^p([0,\w{T}_\ep]\times \R_+), 
\]
where 
\[
\w{T}_\ep=\left\{\begin{array}{ll}
c_0\ep^{-2p(p-1)/\gamma(p,n)}&\mbox{if}\quad \dy\frac{n+1}{n-1}\le p<p_0(n),\\
\exp(c_0\ep^{-p(p-1)})&\mbox{if}\quad p=p_0(n),
\end{array}\right., 
\]
provided $c_0$ is a sufficiently small positive constant. 
Namely, $\w{T}_\ep$ implies the lower bound of the lifespan for (\ref{single}). 
Note that our estimates for $p=q$ in Theorem \ref{main} have the same form as their $\w{T}_\ep$, since we have 
\[
F(p,p;n)=\frac{\gamma(p,n)}{2p(p-1)}.
\]
\end{re}
\par
This article is organized as follows. 
In section \ref{sec:rep}, we recall the representation formulas introduced by Kubo and Kubota in \cite{KK_even}. 
We will consider the decay estimates of a solution to free wave equations in section \ref{sec:free}. 
In section \ref{sec:apriori}, we will state a priori estimates. 
Using them, we will construct a local solution and estimate the lifespan in section \ref{sec:lifespan}. 
Finally section \ref{sec:prop} is devoted to the proof of a priori estimates.

\section{Representation formula}
\label{sec:rep}
\par 
Theorem \ref{main} will be proved by the construction of a local-in-time solution to (\ref{I_rad}) (\ref{ID_rad}) by the iteration argument in a weighted $L^\infty$ space with the sharp pointwise estimates of a solution. 
This is based on the way by Kubo and Kubota \cite{KK_even}. 
In \cite{KK_even}, they get two useful representation formulas starting from the one by Courant and Hilbert \cite{CH}. 
We shall use their formulas and apply their methods into our problem for the compactly supported data. 
\par 
Let $n=2m+2$ $(m\ge 2)$. 
We define in (\ref{I_rad}) that 
\begin{align}
u^0(r,t)&:=\frac{1}{c_n}\left\{\Theta(g_1)(r,t)+D_t\Theta(f_1)(r,t)\right\}, \label{u0}\\
v^0(r,t)&:=\frac{1}{c_n}\left\{\Theta(g_2)(r,t)+D_t\Theta(f_2)(r,t)\right\}, \label{v0}\\
L(G)(r,t)&:=\frac{1}{c_n}\int_0^t\Theta(G(\cdot,\tau))(r,t-\tau)d\tau, \label{L}
\end{align}
where $c_n:=\sqrt{\pi}\Gamma\left((n-1)/2\right)$ and 
\begin{equation}
\label{theta12}
\Theta(g)(r,t):=\frac{1}{r^{2m}}\left(J_1(g)(r,t)+J_2(g)(r,t)\right). 
\end{equation}
Here $J_1$ and $J_2$ are defied by  
\begin{align}
J_1(r,t)&:=\int_{|t-r|}^{t+r}\lambda^{2m+1}g(\lambda)K_m(\lambda,r,t)d\lambda,\label{J1}\\
J_2(r,t)&:=\int_0^{(t-r)_+}\lambda^{2m+1}g(\lambda)\w{K}_m(\lambda,r,t)d\lambda,\label{J2}\\
\end{align}
where
\begin{align}
K_j(\lambda,r,t):=\int_\lambda^{t+r}\frac{H_j(r,t)(\rho,r,t)}{\sqrt{\rho^2-\lambda^2}}d\rho\quad\mbox{for}\quad 0\le j\le m,\label{K}\\
\w{K}_j(\lambda,r,t):=\int_{t-r}^{t+r}\frac{H_j(r,t)(\rho,r,t)}{\sqrt{\rho^2-\lambda^2}}d\rho\quad\mbox{for}\quad 0\le j\le m, \label{Ktilde}
\end{align}
with 
\[
H_j(\rho,r,t):=\left\{D_\rho\left(-\frac{1}{2\rho}\right)\right\}^j(r^2-(\rho-t)^2)^{m-\frac{1}{2}}\qquad\mbox{for}\quad |\rho-t|<r,\ 0\le j\le m. 
\]
This expression was introduced by \cite{KK_even}, while they denote $J_i$ by $W_i$ in \cite{KK_even}. 
First we mention the following property. 
\begin{prop}
Let the assumptions of Theorem \ref{main} be fulfilled. 
Let $(u(r,t),v(r,t))$ be $C^1$ solution of (\ref{I_rad}) with (\ref{ID_rad}). 
If we set $\varphi_i(x)=f_i(|x|)$, $\psi_i(x)=g_i(|x|)$ ($i=1,2$), $\w{u}(x,t)=u(|x|,t)$, $\w{v}(x,t)=v(|x|,t)$ for $(x,t)\in(\R^n\backslash\{0\})\times[0,T]$, 
then, $(\w{u},\w{v})$ is a weak solution of (\ref{D}) with (\ref{ID}) in a distribution sense and it belongs to $((C^2([0,T];\mathcal{D}'(\R^n)))^2$. 
\end{prop}
{\it (Proof)}\\
 This can be proved by the same way as Theorem 2.1 and Theorem 2.2 in \cite{KK_even}. 
So we shall omit the proof. \hfill $\Box$\\
\par
In what follows, we shall summarize the properties of $K_j$ in (\ref{K}) and $\w{K}_j$ in (\ref{Ktilde}) mentioned by Kubo and Kubota in \cite{KK_even} which will be used in our problem. 

\begin{lm}
(Kubo and Kubota \cite{KK_even} \S 4.1 p.215)\\
For $0\le j\le m-1$, we have
\begin{align}
K_j(\lambda,r,t)&=\w{K}(\lambda,r,t)\quad\mbox{for}\quad \lambda=t-r>0,\label{4.41}\\
K_j(\lambda,r,t)&=0\quad\mbox{for}\quad \lambda=t+r.\label{4.42}
\end{align}
\end{lm}

\begin{lm}
(Kubo and Kubota \cite{KK_even} \S 4.1 Lemma 4.1)\\
Let $1\le j\le m$ and $(r,t)\in\Omega$. Then we have
\begin{align}
K_j(\lambda,r,t)&=-D_{\lambda^2}K_{j-1}(\lambda,r,t)\quad\mbox{for}\quad |t-r|<\lambda<t+r,\label{4.1}\\
\w{K}_j(\lambda,r,t)&=-D_{\lambda^2}\w{K}_{j-1}(\lambda,r,t)\quad\mbox{for}\quad 0<\lambda<t-r,\label{4.2}
\end{align}
\end{lm}

\begin{lm}
(Kubo and Kubota \cite{KK_even} \S 4.1 Lemma 4.2)\\
Let $\eta=0\ or\ \frac{1}{2}$ and let $(r,t)\in\Omega$.
Denote $(D_r,D_t,D_\lambda)$ by $D$. 
If $j+|\alpha|\le m$, we have 
\begin{align}
|D^\alpha K_j(\lambda,r,t)
&\le Cr^{2m-j-|\alpha|+\eta-\frac{1}{2}}\lambda^{-j-\eta}(\lambda-t+r)^{-\frac{1}{2}}\quad\mbox{for}\quad |t-r|<\lambda<t+r,\quad t\ge 0,\label{4.4}\\
|D^\alpha \w{K}_j(\lambda,r,t)
&\le Cr^{2m-j-|\alpha|+\eta-\frac{1}{2}}(t-r)^{-j-\eta}(t-r-\lambda)^{-\frac{1}{2}}\quad\mbox{for}\quad 0<\lambda<t-r.\label{4.5}
\end{align}
\end{lm}

\begin{lm}
(Kubo and Kubota \cite{KK_even} \S 4.1 Corollary 4.3)\\
Let $\eta=0\ or\ \frac{1}{2}$ and let $(r,t)\in\Omega$. 
If $j+|\alpha|\le m$, we have 
\begin{equation}
|D_\lambda D_{r,t}^\alpha \w{K}_j(\lambda,r,t)|
\le Cr^{2m-j-|\alpha|+\eta-\frac{1}{2}}(t-r)^{-j-\eta}(t-r-\lambda)^{-\frac{3}{2}}\quad\mbox{for}\quad 0<\lambda<t-r.
\label{4.13}
\end{equation}
\end{lm}

\begin{lm}
(Kubo and Kubota \cite{KK_even} \S 4.1 Lemma 4.4)\\
Let $(r,t)\in\Omega$ and $0<\lambda<t-r$. 
Let $\eta=0\ or\ \frac{1}{2}$. 
If $0\le j\le m$, $|\alpha|\le m$ and $j+|\alpha|\ge m$, we have 
\begin{equation}
|D_{r,t}^\alpha \w{K}_j(\lambda,r,t)|
\le Cr^{m-\eta}(t-r)^{-j-\frac{1}{2}+\eta}(t-r-\lambda)^{m-j-|\alpha|-\frac{1}{2}}. 
\label{4.14}
\end{equation}
Moreover, if $0\le j\le m$, $|\alpha|\le m-1$ and $j+|\alpha|\ge m-1$, we have 
\begin{equation}
|D_{r,t}^\alpha \w{K}_j(\lambda,r,t)|
\le Cr^{m+1-\eta}(t-r)^{-j-\frac{1}{2}+\eta}(t-r-\lambda)^{m-j-|\alpha|-\frac{3}{2}}.
\label{4.15}
\end{equation}
\end{lm}

\begin{lm}
(Kubo and Kubota \cite{KK_even} \S 4.1 Lemma 4.5)\\
Let $\eta=0\ or\ \frac{1}{2}$ and let $(r,t)\in\Omega$. 
If $0<r-t\le\lambda<r+t$, we have
\begin{equation}
|D_{r,t,\lambda}^\alpha K_m(\lambda,r,t)|
\le Cr^{m-\frac{1}{2}+\eta}\lambda^{-m-|\alpha|-\frac{1}{2}-\eta.}
\label{4.18}
\end{equation}
\end{lm}

\begin{lm}
(Kubo and Kubota \cite{KK_even} \S 4.1 p.213 REMARK)\\
When $0<t<r$, we have for $\lambda=r+t$ 
\begin{align}
|K_m(\lambda,r,t)|&\le Cr^{m-\frac{1}{2}+\mu}\lambda^{-m-\frac{1}{2}-\mu}\label{4.30}\\
|D_{r,t,\lambda}K_m(\lambda,r,t)|&\le Cr^{m-\frac{1}{2}}\lambda^{-m-\frac{3}{2}}. \label{4.32}
\end{align}
\end{lm}

\par
Using these properties, Kubo and Kubota obtained another representation formula for $\Theta$ as follows by Lemma 4.6 and Lemma 4.7 in \cite{KK_even}.  
\par
First, assume that 
\begin{equation}
\label{A_g}
g\in C^1((0,\infty))\quad\mbox{and}\quad 
g^{(j)}(\lambda)=O(\lambda^{-2m+\delta-j})\quad\mbox{as}\quad \lambda\downarrow 0
\end{equation}
hold for $j=0,1$ and some positive constant $\delta$. 
Then we have 
\begin{equation}
\label{theta34}
2r^{2m}\Theta(g)(r,t)=J_3(r,t)+J_4(r,t),
\end{equation}
where 
\begin{align}
J_3(r,t)&:=\int_{|t-r|}^{t+r}(\lambda^{2m}g(\lambda))'K_{m-1}(\lambda,r,t),\label{J3}\\
J_4(r,t)&:=\dy\left\{\begin{array}{ll}
\dy\int_0^{t-r}(\lambda^{2m}g(\lambda))'\w{K}_{m-1}(\lambda,r,t)d\lambda,&\mbox{if}\quad t-r>0,\\
(r-t)^{2m}g(r-t)K_{m-1}(r-t,r,t)&\mbox{if}\quad t-r<0.
\end{array}\right.\label{J4}
\end{align}
Moreover, we have 
\begin{equation}
\label{theta56}
D_\alpha(2r^{2m}\Theta(g)(r,t))=J_5(r,t)+J_6(r,t),
\end{equation}
where $D_\alpha=D_r\ or\ D_t$ and 
\begin{align}
J_5(r,t)&:=\int_{|t-r|}^{t+r}(\lambda^{2m}g(\lambda))'D_\alpha K_{m-1}(\lambda,r,t),\label{J5}\\
J_6(r,t)&:=\dy\left\{\begin{array}{ll}
\dy\int_0^{t-r}(\lambda^{2m}g(\lambda))'D_\alpha\w{K}_{m-1}(\lambda,r,t)d\lambda,&\mbox{if}\quad t-r>0,\\
(r-t)^{2m}g(r-t)D_\alpha K_{m-1}(r-t,r,t)&\mbox{if}\quad t-r<0.
\end{array}\right.\label{J6}
\end{align}
\par
Next, assume that 
\begin{equation}
\label{A_f}
f\in C^2((0,\infty))\quad\mbox{and}\quad 
f^{(j)}(\lambda)=O(\lambda^{-2m+2+\delta-j})\quad\mbox{as}\quad \lambda\downarrow 0
\end{equation}
for $j=0,1,2$ and some positive constant $\delta$. 
When $t-r>0$, we have 
\begin{equation}
\label{theta78}
D^\beta D_t(2r^{2m}\Theta(f)(r,t))=J_7(r,t)+J_8(r,t),
\end{equation}
where $D=(D_r,D_t)$ and 
\begin{align}
J_7(r,t)&:=\int_{t-r}^{t+r}D_\lambda D_{\lambda^2}(\lambda^{2m}f(\lambda))\cdot D^\beta D_tK_{m-2}(\lambda,r,t),\label{J7}\\
J_8(r,t)&:=\int_0^{t-r}D_\lambda D_{\lambda^2}(\lambda^{2m}f(\lambda))\cdot D^\beta D_t\w{K}_{m-2}(\lambda,r,t).\label{J8}\\
\end{align}
On the other hand, when $t-r<0$, we have 
\begin{equation}
\label{theta910}
D_t(2r^{2m}\Theta(f)(r,t))=J_9(r,t)+J_{10}(r,t),
\end{equation}
where
\begin{align}
J_9(r,t)&:=\int_{r-t}^{r+t}\lambda^{2m+1}f(\lambda)\cdot D_tK_{m}(\lambda,r,t),\label{J9}\\
J_{10}(r,t)&:=(\lambda^{2m+1}f(\lambda)K_m(\lambda,r,t))|_{\lambda=r+t}+(\lambda^{2m+1}f(\lambda)K_m(\lambda,r,t))|_{\lambda=r-t}.\label{J10}
\end{align}
\par
Our assumptions on $g_i$ and $f_i$ ($i=1,2$) in Theorem \ref{main} 
can play the same roles as (\ref{A_g}) and (\ref{A_f}) in the proof of the above representation formulas in 
Kubo and Kubota \cite{KK_even}. 
Therefore, we will use these representation formulas (\ref{theta12}), (\ref{theta34}), (\ref{theta56}), (\ref{theta78}), (\ref{theta910}) suitably in the forthcoming sections.

\section{Decay estimates of free solution}

\label{sec:free}
\par
In this section, we show the decay estimates of $u^0$ and $v^0$. 
We have the following. 

\begin{prop}
\label{prop:free}
Let the assumptions on $f_i$, $g_i$ ($i=1,2$) in Theorem \ref{main} be fulfilled. 
Let $u^0$ and $v_0$ are the one in (\ref{u0}) and (\ref{v0}), respectively. 
Then we have $\supp\{u^0,v^0\}\subset\{r\le t+k\}$ and 
\begin{align*}
|D^\beta u^0(r,t)|&\le \w{C}\left(\frac{r}{k}\right)^{-m+1-|\beta|}\left(\frac{r+2k}{k}\right)^{-1+|\beta|}\tau_+(r,t)^{-\frac{1}{2}}\phi_i(r,t), \\
|D^\beta v^0(r,t)|&\le \w{C}\left(\frac{r}{k}\right)^{-m+1-|\beta|}\left(\frac{r+2k}{k}\right)^{-1+|\beta|}\tau_+(r,t)^{-\frac{1}{2}}\phi_i(r,t), 
\end{align*}
for $(r,t)\in\Omega:=\{(r,t)\in(0,\infty)\times[0,\infty) : r\le t+k\}$, 
where $D=(D_r,D_t)$, $|\beta|\le 1$, $\tau_{\pm}(r,t)=(t\pm r+2k)/k$, $\phi_i(r,t)=\tau_-^{-m-\frac{1}{2}-|\beta|}$ and $\w{C}$ is a positive constant independent of $\ep$.
\end{prop}
{\bf Proof.} 
We shall consider $u^0$ only, because $v_0$ has the same form. 
The support property of $u^0$ can be readily obtained by the representation formulra (\ref{u0}), (\ref{theta12}), (\ref{J1}) and (\ref{J2}). 
In what follows, we prove the decay estimates of $D^\beta u^0$ by the dividng into three cases. 
This division is due Kubo and Kubota in (\cite{KK_even}). 
Here and hereafter, we denote any positive constant independent of $\ep$ by $C$. 
Namely, $C$ can mean the different values in each line. 
Besides, we denote $\tau_{\pm}(r,t)$ by $\tau_{\pm}$ simply, since we fix $(r,t)$ in each cases. 
Similarly, we shall omit the variable $(r,t)$ of the other functions as well when not misleading. 
\subsection{Case 1, $t\ge 2r>0$ or $0<r<\min\{t,k\}$}
In this case, it holds that 
\begin{align*}
t\ge 2r>0 &\Longrightarrow 3\tau_-\ge\tau_+\ge\frac{r+2k}{k},\\
0<r<\min\{t,r\}&\Longrightarrow 2<\frac{r+2k}{k},\ \tau_+,\ \tau_-\le 5, 
\end{align*}
for $(r,t)\in\Omega$. 
Hence, it is enough to show 
\begin{align}
|D^\beta D_t(2r^{2m}\Theta(f)(r,t))|
&\le Cr^{m+1-|\beta|}\tau_-^{-m-2},\label{decay1_f}\\
|D^\beta (2r^{2m}\Theta(g)(r,t))|
&\le Cr^{m+1-|\beta|}\tau_-^{-m-2},\label{decay1_g}
\end{align}
for $|\beta|=0,1$ in $\Omega$. 
\par
For (\ref{decay1_f}), we use the representation formula (\ref{theta78}) both for $|\beta|=0$ and $1$. 
By (\ref{4.4}) with $\eta=0$, we have 
\[
|J_7(r,t)|
\le Cr^{m+\frac{1}{2}-|\beta|}\int_{t-r}^{t+r}\frac{\lambda^{m-1}|f(\lambda)|+\lambda^{m}|f'(\lambda)|+\lambda^{m+1}|f''(\lambda)|}{\sqrt{\lambda-t+r}}d\lambda.
\]
When $t-r\ge k$, $J_7=0$ because of the assumption of the support property of data, {\it i.e.}
\begin{equation}
\label{s}
f(r)=g(r)=0\qquad\mbox{if}\quad r\ge k.
\end{equation}
When $0<t-r\le k$, it follows from (\ref{s}) and $2\le\tau_-\le 3$ that 
\begin{align*}
|J_7(r,t)|
&\le Cr^{m+\frac{1}{2}-|\beta|}\int_{t-r}^{\min\{k,t+r\}}\frac{1}{\sqrt{\lambda-t+r}}d\lambda
\le Cr^{m+1-|\beta|}\tau_-^{-m-2}.
\end{align*} 
As to $J_8$ for $0<t-r\le 2k$, we use $(\ref{4.5})_{m-2}$ with $\mu=\frac{1}{2}$. 
Then, it follows from (\ref{s}) and $2\le\tau_-\le 4$ that 
\begin{align*}
|J_8(r,t)|
&\le Cr^{m+1-|\beta|}(t-r)^{-m+\frac{3}{2}}\int_0^{t-r}\frac{\lambda^{2m-3}|f(\lambda)|+\lambda^{2m-2}|f'(\lambda)|+\lambda^{2m-1}|f''(\lambda)|}{\sqrt{t-r-\lambda}}d\lambda\\
&\le Cr^{m+1-|\beta|}\int_0^{t-r}\frac{\lambda^{m-\frac{3}{2}}|f(\lambda)|+\lambda^{m-\frac{1}{2}}|f'(\lambda)|+\lambda^{m+\frac{1}{2}}|f''(\lambda)|}{\sqrt{t-r-\lambda}}d\lambda\\
&\le Cr^{m+1-|\beta|}\int_0^{t-r}\frac{1}{\sqrt{t-r-\lambda}}d\lambda\\
&\le Cr^{m+1-|\beta|}\tau_-^{-m-2}.
\end{align*}
When $t-r\ge 2k$, (\ref{s}) gives us 
\begin{align*}
J_8(r,t)
&=\int_0^{\frac{t-r}{2}}D_\lambda D_{\lambda^2}(\lambda^{2m}f(\lambda))\cdot D^\beta D_t\w{K}_{m-2}(\lambda,r,t)d\lambda. 
\end{align*}
By employing the integration by parts and (\ref{4.2}) twice together with (\ref{s}), we have 
\begin{align}
J_8(r,t)
&=\left[D_{\lambda^2}(\lambda^{2m}f(\lambda))D^\beta D_t\w{K}_{m-2}(\lambda,r,t)\right]_{\lambda=0}^{\frac{t-r}{2}}\nonumber\\
&\quad+\left[\lambda^{2m}f(\lambda) D^\beta D_t\w{K}_{m-1}(\lambda,r,t)\right]_{\lambda=0}^{\frac{t-r}{2}}\nonumber\\
&\quad+2\int_0^{\frac{t-r}{2}}\lambda^{2m+1}f(\lambda)D^\beta D_t\w{K}_{m}(\lambda,r,t)d\lambda.\label{parts}
\end{align}
By (\ref{s}), the first and the second line of (\ref{parts}) vanishes. 
Using $(\ref{4.15})_m$ with $\mu=0$ for $|\beta|=0$ and $(\ref{4.14})_m$ with $\mu=0$ for $|\beta|=1$, we have 
\begin{align*}
|J_8(r,t)|
&\le Cr^{m+1-|\beta|}(t-r)^{-m-3}\int_0^{\frac{t-r}{2}}\lambda^{2m+1}|f(\lambda)|d\lambda
\le Cr^{m+1-|\beta|}\tau_-^{-m-2}.
\end{align*}
\par
For (\ref{decay1_g}) with $|\beta|=0$, we use (\ref{theta34}). 
Using $(\ref{4.4})_{m-1}$ with $\mu=0$, we have 
\[
|J_3(r,t)|
\le Cr^{m+\frac{1}{2}}\int_{t-r}^{t+r}\frac{\lambda^m|g(\lambda)|+\lambda^{m+1}|g'(\lambda)|}{\sqrt{\lambda-t+r}}d\lambda. 
\]
Then, we can treat it by the same way as $J_7$ for (\ref{decay1_f}) and get the desired estimate. 
As to $J_4$ for $t-r\le 2k$, we use $(\ref{4.5})_{m-1}$ with $\mu=\frac{1}{2}$. 
Then, it follows from $2\le\tau_-\le 4$ that 
\begin{align*}
|J_4(r,t)|
&\le Cr^{m+1}(t-r)^{-m+\frac{1}{2}}\int_0^{t-r}\frac{\lambda^{2m+1}|g(\lambda)|+\lambda^{2m}|g'(\lambda)|}{\sqrt{t-r-\lambda}}d\lambda\\
&\le Cr^{m+1}(t-r)^{\frac{1}{2}}\int_0^{t-r}\frac{\lambda^{m+1}|g(\lambda)|+\lambda^{m}|g'(\lambda)|}{\sqrt{t-r-\lambda}}d\lambda\\
&\le Cr^{m+1}(t-r)^{\frac{1}{2}}\int_0^{t-r}\frac{1}{\sqrt{t-r-\lambda}}d\lambda\\
&\le Cr^{m+1}\tau_-^{-m-2}.
\end{align*}
When $t-r\ge 2k$, using (\ref{s}), the integration by parts and (\ref{4.2}), we have 
\begin{align*}
J_4(r,t)
&=\int_0^{\frac{t-r}{2}}(\lambda^{2m}g(\lambda))'\w{K}_{m-1}(\lambda,r,t)d\lambda\\
&=-\int_0^{\frac{t-r}{2}}\lambda^{2m}g(\lambda)D_\lambda\w{K}_{m-1}(\lambda,r,t)d\lambda\\
&=-2\int_0^{\frac{t-r}{2}}\lambda^{2m+1}g(\lambda)D_{\lambda^2}\w{K}_{m-1}(\lambda,r,t)d\lambda\\
&=2\int_0^{\frac{t-r}{2}}\lambda^{2m+1}g(\lambda)\w{K}_{m}(\lambda,r,t)d\lambda. 
\end{align*}
By$(\ref{4.15})_{m}$ with $\mu=0$, we get 
\begin{align*}
|J_4(r,t)|
&\le Cr^{m+1}(t-r)^{-m-\frac{1}{2}}\int_0^{\frac{t-r}{2}}\frac{\lambda^{2m+1}|g(\lambda)|}{(t-r-\lambda)^{\frac{3}{2}}}d\lambda\\
&\le Cr^{m+1}(t-r)^{-m-2}\int_0^k\lambda^{2m+1}|g(\lambda)|d\lambda\\
&\le Cr^{m+1}\tau_-^{-m-2}.
\end{align*}
\par
For (\ref{decay1_g}) with $|\beta|=1$, we use (\ref{theta56}). 
If we use
\begin{center}
\begin{tabular}{rl}
$(\ref{4.4})_{m-1}$ with $\mu=0$ &for $J_5$, \\
$(\ref{4.5})_{m-1}$ with $\mu=\frac{1}{2}$ &for $J_6$ with $t-r\le 2k$,\\
(\ref{4.2}) and $(\ref{4.14})_m$ with $\mu=0$ &for $J_6$ with $t-r\ge 2k$,
\end{tabular}
\end{center}
we can treat $J_5$ and $J_6$ by the same way as $J_7$ and $J_8$ for (\ref{decay1_f}), respectively. 

\subsection{Case 2, $k\le r\le t\le 2r$}
\par
In this case, it holds that 
\[
\tau_-\le \tau_+\le 3\frac{r+2k}{k}\le\frac{9r}{k}\\
\]
for $(r,t)\in\Omega$. 
Hence, it is enough to show 
\begin{align}
|D^\beta D_t(2r^{2m}\Theta(f)(r,t))|
&\le Cr^{m-\frac{1}{2}}\tau_-^{-m-\frac{1}{2}-|\beta|},\label{decay2_f}\\
|D^\beta (2r^{2m}\Theta(g)(r,t))|
&\le Cr^{m-\frac{1}{2}}\tau_-^{-m-\frac{1}{2}-|\beta|},\label{decay2_g}
\end{align}
for $|\beta|=0,1$ in $\Omega$. 
\par
For (\ref{decay2_f}) with $|\beta|=0$, we use (\ref{theta56}). 
Using $(\ref{4.4})_{m-1}$ with $\mu=0$, 
\begin{align*}
J_5(r,t)
&= Cr^{m-\frac{1}{2}}\int_{t-r}^{t+r}
\frac{\lambda^{m}|f(\lambda)|+\lambda^{m+1}|f'(\lambda)|}{\sqrt{\lambda-t+r}}d\lambda.\\
\end{align*}
Since $J_5(r,t)=0$ for $t-r\ge k$ by (\ref{s}), we consider it when $t-r\le k$. 
Then, we have by (\ref{s})
\begin{align*}
J_5(r,t)
&\le Cr^{m-\frac{1}{2}}\\
&\le Cr^{m-\frac{1}{2}}\tau_-^{-m-\frac{1}{2}}. 
\end{align*}
Next we consider $J_6$. 
When $t-r\le 2k$, using $(\ref{4.14})_{m-1}$ with $\mu=\frac{1}{2}$ and (\ref{s}), we have 
\begin{align*}
J_6(r,t)
&\le Cr^{m-\frac{1}{2}}(t-r)^{-m+1}\int_0^{t-r}\frac{\lambda^{2m-1}|f(\lambda)|+\lambda^{2m}|f'(\lambda)|}{\sqrt{t-r-\lambda}}d\lambda\\
&\le Cr^{m-\frac{1}{2}}(t-r)\int_0^{t-r}\frac{\lambda^{m-1}|f(\lambda)|+\lambda^{m}|f'(\lambda)|}{\sqrt{t-r-\lambda}}d\lambda\\
&\le Cr^{m-\frac{1}{2}}\\
&\le Cr^{m-\frac{1}{2}}\tau_-^{-m-\frac{1}{2}}. 
\end{align*}
When $t-r\ge 2k$, by the integration by parts, (\ref{s}) and (\ref{4.2}), we have 
\begin{align*}
J_6(r,t)
&=\int_0^{\frac{t-r}{2}}(\lambda^{2m}f(\lambda))'D_t\w{K}_{m-1}(\lambda,r,t)d\lambda\\
&=-\int_0^{\frac{t-r}{2}}\lambda^{2m}f(\lambda)D_\lambda D_t\w{K}_{m-1}(\lambda,r,t)d\lambda.\\
&= 2\int_0^{\frac{t-r}{2}}\lambda^{2m+1}f(\lambda)D_t\w{K}_{m}(\lambda,r,t)d\lambda.
\end{align*}
Hence, by $(\ref{4.14})_{m-1}$ with $\eta=\frac{1}{2}$, we have  
\begin{align*}
|J_6(r,t)|
&\le Cr^{m-\frac{1}{2}}(t-r)^{-m}\int_0^{\frac{t-r}{2}}\frac{\lambda^{2m+1}f(\lambda)}{(t-r-\lambda)^{\frac{3}{2}}}d\lambda\\
&\le Cr^{m-\frac{1}{2}}(t-r)^{-m-\frac{3}{2}}\int_0^{\frac{t-r}{2}}\lambda^{2m+1}f(\lambda)d\lambda\\
&\le Cr^{m-\frac{1}{2}}\tau_-^{-m-\frac{1}{2}}. 
\end{align*}
\par
For (\ref{decay2_f}) with $|\beta|=1$, we use (\ref{theta78}). 
Using $(\ref{4.4})_{m-2}$ with $\mu=0$, we have 
\begin{align*}
|J_7(r,t)|
&\le Cr^{m-\frac{1}{2}}\int_{t-r}^{t+r}
\frac{\lambda^{m-1}|f(\lambda)|+\lambda^{m}|f'(\lambda)|+\lambda^{m+1}|f''(\lambda)|}
{\sqrt{\lambda-t+r}}d\lambda.
\end{align*}
When $t-r\ge k$, $J_7=0$. 
When $t-r\le k$, we can treat it as before. 
\par
As to $J_8$, dividing it as 
\[
\int_0^{t-r}d\lambda=\int_0^{\frac{t-r}{2}}d\lambda+\int_{\frac{t-r}{2}}^{t-r}d\lambda. 
\]
and calculating the first term by the same way as in (\ref{parts}), 
we have 
\begin{align*}
J_8(r,t)
&=\left[D_{\lambda^2}(\lambda^{2m}f(\lambda))D^\beta D_t\w{K}_{m-2}(\lambda,r,t)\right]_{\lambda=0}^{\frac{t-r}{2}}\nonumber\\
&\quad+\left[\lambda^{2m}f(\lambda) D^\beta D_t\w{K}_{m-1}(\lambda,r,t)\right]_{\lambda=0}^{\frac{t-r}{2}}\nonumber\\
&\quad+2\int_0^{\frac{t-r}{2}}\lambda^{2m+1}f(\lambda)D^\beta D_t\w{K}_{m}(\lambda,r,t)d\lambda\\
&\quad+\int_{\frac{t-r}{2}}^{t-r}D_\lambda D_{\lambda^2}(\lambda^{2m}f(\lambda))\cdot D^\beta D_t \w{K}_{m-2}(\lambda,r,t)d\lambda\\
&\equiv A_1+A_2+A_3+A_4. 
\end{align*}
For $A_1$, using $(\ref{4.5})_{m-2}$ with $\eta=0$ to the term of $\lambda=\frac{t-r}{2}$, we have 
\begin{align*}
|A_1|
&\le Cr^{m-\frac{1}{2}}\left((t-r)^{m-\frac{1}{2}}\left|f\left(\frac{t-r}{2}\right)\right|+(t-r)^{m+\frac{1}{2}}\left|f'\left(\frac{t-r}{2}\right)\right|\right).
\end{align*}
When $t-r\ge 2k$, $A_1=0$. 
When $t-r\le 2k$, we have $|A_1|\le Cr^{m-\frac{1}{2}}$. \\
For $A_2$, using (\ref{4.14}) with $\eta=\frac{1}{2}$ to the term of $\lambda=\frac{t-r}{2}$, we have 
\begin{align*}
|A_2|
&\le Cr^{m-\frac{1}{2}}(t-r)^{m-\frac{1}{2}}\left|f\left(\frac{t-r}{2}\right)\right|.
\end{align*}
When $t-r\ge 2k$, $A_2=0$. 
When $t-r\le 2k$, we have 
$|A_2|\le Cr^{m-\frac{1}{2}}$\\
For $A_3$, using (\ref{4.14}) with $\eta=\frac{1}{2}$, we have 
\begin{align*}
|A_3|
&\le Cr^{m-\frac{1}{2}}(t-r)^{-m}\int_0^{\frac{t-r}{2}}\frac{\lambda^{2m+1}|f(\lambda)|}{(t-r-\lambda)^\frac{5}{2}}d\lambda\\
&\le Cr^{m-\frac{1}{2}}(t-r)^{-m-\frac{5}{2}}\int_0^{\frac{t-r}{2}}\lambda^{2m+1}|f(\lambda)|d\lambda.
\end{align*}
When $t-r\ge 2k$, we have 
\begin{align*}
|A_3|
&\le Cr^{m-\frac{1}{2}}(t-r)^{-m-\frac{5}{2}}\int_0^k\lambda^{2m+1}|f(\lambda)|d\lambda\\
&\le Cr^{m-\frac{1}{2}}\tau_-^{-m-\frac{5}{2}}.
\end{align*}
When $t-r\le 2k$, we have 
\begin{align*}
|A_3|
&\le Cr^{m-\frac{1}{2}}(t-r)^{m-\frac{3}{2}}\int_0^{\frac{t-r}{2}}|f(\lambda)|d\lambda\\
&\le Cr^{m-\frac{1}{2}}.
\end{align*}
For $A_4$, using $(\ref{4.5})_{m-2}$ with $\eta=0$, we have 
\begin{align*}
A_4
&\le Cr^{m-\frac{1}{2}}(t-r)^{-m+2}\int_{\frac{t-r}{2}}^{t-r}
\frac{\lambda^{2m-3}|f(\lambda)|+\lambda^{2m-2}|f'(\lambda)|+\lambda^{2m-1}|f''(\lambda)|}
{\sqrt{t-r-\lambda}}d\lambda. 
\end{align*}
When $t-r\ge 2k$, $B=0$. 
When $t-r\le 2k$, we have 
\begin{align*}
A_4
&\le Cr^{m-\frac{1}{2}}\int_{\frac{t-r}{2}}^{t-r}
\frac{\lambda^{m-1}|f(\lambda)|+\lambda^{m}|f'(\lambda)|+\lambda^{m+1}|f''(\lambda)|}
{\sqrt{t-r-\lambda}}d\lambda\\
&\le Cr^{m-\frac{1}{2}}\int_{\frac{t-r}{2}}^{t-r}
\frac{1}{\sqrt{t-r-\lambda}}d\lambda\\
&\le Cr^{m-\frac{1}{2}}
\end{align*}
Therefore we get 
\begin{align*}
|J_8(r,t)|
&\le\left\{\begin{array}{ll}
Cr^{m-\frac{1}{2}}\tau_-^{-m-\frac{5}{2}}&\mbox{if}\quad t-r\ge 2k,\\
Cr^{m-\frac{1}{2}}&\mbox{if}\quad t-r\le 2k,
\end{array}\right.\\
&\le Cr^{m-\frac{1}{2}}\tau_-^{-m-\frac{3}{2}}. 
\end{align*}
\par
For (\ref{decay2_g}) with $|\beta|=0$, we use (\ref{theta12}). 
Usin $(\ref{4.4})_m$ with $\eta=0$, we have 
\[
|J_1(r,t)|\le Cr^{m-\frac{1}{2}}\int_{t-r}^{t+r}\frac{\lambda^{m+1}|g(\lambda)|}{\sqrt{\lambda-t+r}}d\lambda. 
\]
When $t-r\ge k$, $J_1=0$. 
When $t-r\le k$, we have 
\begin{align*}
|J_1(r,t)|
&\le Cr^{m-\frac{1}{2}}\int_{t-r}^k\frac{1}{\sqrt{\lambda-t+r}}d\lambda\\
&\le Cr^{m-\frac{1}{2}}\\
&\le Cr^{m-\frac{1}{2}}\tau_-^{-m-\frac{1}{2}}.
\end{align*}
Using $(\ref{4.5})_m$ with $\eta=0$, we have 
\[
|J_2(r,t)|
\le Cr^{m-\frac{1}{2}}(t-r)^{-m}
\int_0^{t-r}\frac{\lambda^{2m+1}|g(\lambda)|}{\sqrt{t-r-\lambda}}d\lambda. 
\]
When $t-r\le 2k$, 
\begin{align*}
|J_2(r,t)|
&\le Cr^{m-\frac{1}{2}}\int_0^{t-r}\frac{\lambda^{m+1}|g(\lambda)|}{\sqrt{t-r-\lambda}}d\lambda\\
&\le Cr^{m-\frac{1}{2}}\int_0^{t-r}\frac{1}{\sqrt{t-r-\lambda}}d\lambda\\
&\le Cr^{m-\frac{1}{2}}\\
&\le Cr^{m-\frac{1}{2}}\tau_-^{-m-\frac{1}{2}}. 
\end{align*}
When $t-r\ge 2k$, we have by (\ref{s})
\begin{align*}
|J_2(r,t)|
&\le Cr^{m-\frac{1}{2}}(t-r)^{-m}\int_0^{\frac{t-r}{2}}\frac{\lambda^{m+1}|g(\lambda)|}{\sqrt{t-r-\lambda}}d\lambda\\
&\le Cr^{m-\frac{1}{2}}(t-r)^{-m-\frac{1}{2}}\int_0^k\lambda^{m+1}|g(\lambda)|d\lambda\\
&\le Cr^{m-\frac{1}{2}}\tau_-^{-m-\frac{1}{2}}. 
\end{align*}
\par
For (\ref{decay2_g}) with $|\beta|=1$, we use (\ref{theta56}). 
Using $(\ref{4.4})_{m-1}$ with $\eta=0$, we have 
\[
|J_5(r,t)|
\le Cr^{m-\frac{1}{2}}\int_{t-r}^{t+r}\frac{\lambda^m|g(\lambda)|+\lambda^{m+1}|g'(\lambda)|}{\sqrt{\lambda-t+r}}d\lambda. 
\]
When $t-r\ge k$, $J_5=0$. 
When $t-r\le k$, 
\begin{align*}
|J_5(r,t)|
&\le Cr^{m-\frac{1}{2}}\int_{t-r}^k\frac{1}{\sqrt{\lambda-t+r}}d\lambda\\
&\le Cr^{m-\frac{1}{2}}\\
&\le Cr^{m-\frac{1}{2}}\tau_-^{-m-\frac{3}{2}}.
\end{align*}
\par
As to $J_6$ for $t-r\le 2k$, $(\ref{4.5})_{m-1}$ with $\eta=0$ gives us 
\begin{align*}
|J_6(r,t)|
&\le Cr^{m-\frac{1}{2}}(t-r)^{-m+1}\int_0^{t-r}\frac{\lambda^{2m-1}|g(\lambda)|+\lambda^{2m}|g'(\lambda)|}{\sqrt{t-r-\lambda}}d\lambda\\
&\le Cr^{m-\frac{1}{2}}\int_0^{t-r}\frac{\lambda^{m}|g(\lambda)|+\lambda^{m+1}|g'(\lambda)|}{\sqrt{t-r-\lambda}}d\lambda\\
&\le Cr^{m-\frac{1}{2}}\int_0^{t-r}\frac{1}{\sqrt{t-r-\lambda}}d\lambda\\
&\le Cr^{m-\frac{1}{2}}\\
&\le Cr^{m-\frac{1}{2}}\tau_-^{-m-\frac{3}{2}}.
\end{align*}
When $t-r\ge 2k$, using (\ref{s}), (\ref{4.2}) and the integration by parts, we have 
\begin{align*}
J_6(r,t)
&=\int_0^{\frac{t-r}{2}}(\lambda^{2m}g(\lambda))'D_\alpha\w{K}_{m-1}(\lambda,r,t)d\lambda\\
&=-\int_0^{\frac{t-r}{2}}\lambda^{2m}g(\lambda)D_\lambda D_\alpha\w{K}_{m-1}(\lambda,r,t)d\lambda\\
&=2\int_0^{\frac{t-r}{2}}\lambda^{2m+1}g(\lambda)D_\alpha\w{K}_{m}(\lambda,r,t)d\lambda. 
\end{align*}
Then, by (\ref{4.14}) with $\eta=\frac{1}{2}$, we have 
\begin{align*}
|J_6(r,t)|
&\le Cr^{m-\frac{1}{2}}(t-r)^{-m}\int_0^{\frac{t-r}{2}}\frac{\lambda^{2m+1}|g(\lambda)|}{(t-r-\lambda)^{\frac{3}{2}}}d\lambda\\
&\le Cr^{m-\frac{1}{2}}(t-r)^{-m-\frac{3}{2}}\int_0^{k}\lambda^{2m+1}|g(\lambda)|d\lambda\\
&\le Cr^{m-\frac{1}{2}}\tau_-^{-m-\frac{3}{2}}.
\end{align*}
\subsection{Case 3, $r>t$}
\par
In this case, it holds that $1<\tau_-\le 2$ and 
\begin{align*}
1\le \frac{r+2k}{k},\ \tau_-,\ \tau_+\le 4\qquad&\mbox{if}\quad 0<r\le k,\\
\tau_+\le 2\frac{r+2k}{k}\le \frac{6r}{k}\qquad&\mbox{if}\quad r\ge k,
\end{align*}
for $(r,t)\in\Omega$. 
Hence, it is enough to show 
\begin{align}
|D^\beta D_t(2r^{2m}\Theta(f)(r,t))|
&\le \left\{\begin{array}{ll}
Cr^{m+1-|\beta|}&\mbox{if}\quad 0<r\le k,\\
Cr^{m-\frac{1}{2}}&\mbox{if}\quad r\ge k,
\end{array}\right.\label{decay3_f}\\
|D^\beta (2r^{2m}\Theta(g)(r,t))|
&\le \left\{\begin{array}{ll}
Cr^{m+1-|\beta|}&\mbox{if}\quad 0<r\le k,\\
Cr^{m-\frac{1}{2}}&\mbox{if}\quad r\ge k,
\end{array}\right.\label{decay3_g}
\end{align}
for $|\beta|=0,1$ in $\Omega$. 

\par
For (\ref{decay3_f}), we use (\ref{theta910}) for $|\beta|=0$ and 
and the derivative of (\ref{theta910}) for $|\beta|=1$. 
Using (\ref{4.18}) with $\eta=0$, we have 
\[
|J_9(r,t)|\le Cr^{m-\frac{1}{2}}\int_{r-t}^{r+t}\lambda^{m-\frac{1}{2}}|f(\lambda)|d\lambda. 
\]
When $0<r\le k$, we have 
\[
|J_9(r,t)|
\le Cr^{m-\frac{1}{2}}(r+t)^{\frac{3}{2}}\int_{r-t}^{r+t}\lambda^{m-2}|f(\lambda)|d\lambda
\le Cr^{m+1}. 
\]
When $r\ge k$, we have 
\[
|J_9(r,t)|
\le Cr^{m-\frac{1}{2}}\int_0^k\lambda^{m-\frac{1}{2}}|f(\lambda)|d\lambda
\le Cr^{m-\frac{1}{2}}. 
\]
Next, we consider $D_\alpha J_9(r,t)$. 
Differentiating $J_9$ and using (\ref{4.18}) with $\mu=0$, we have 
\begin{align*}
|D J_9(r,t)|
&\le\int_{r-t}^{r+t}\lambda^{2m+1}|f(\lambda)||(D_\alpha\de_tK_m)(\lambda,r,t)|d\lambda\\
&\quad +(r+t)^{2m+1}|f(r+t)||(\de_t K_m)(r+t,r,t)|\\
&\quad +(r-t)^{2m+1}|f(r-t)||(\de_t K_m)(r-t,r,t)|\\
&\le Cr^{m-\frac{1}{2}}\left\{\int_{r-t}^{r+t}\lambda^{m-\frac{3}{2}}|f(\lambda)|d\lambda
+(r+t)^{m-\frac{1}{2}}|f(r+t)|
+(r-t)^{m-\frac{1}{2}}|f(r-t)|\right\}.
\end{align*}
It follows from (\ref{s}) that 
\begin{align*}
|D J_9(r,t)|
&\le Cr^m\left\{
\int_{r-t}^{r+t}\lambda^{m-2}|f(\lambda)|d\lambda
+(r+t)^{m-1}|f(r+t)|
+(r-t)^{m-1}|f(r-t)|
\right\}\le Cr^m
\end{align*}
for $0<r\le k$ and 
\[
|D J_9(r,t)|\le Cr^{m-\frac{1}{2}}. 
\]
for $r\ge k$. 
\par
Third, we consider $J_{10}$. 
Using (\ref{4.30}) for $K_m(r+t,r,t)$ and (\ref{4.18}) with $\eta=0$ for $K_m(r-t,r,t)$, we have 
\[
|J_{10}(r,t)|
\le Cr^{m-\frac{1}{2}}\{(r+t)^{m+\frac{1}{2}}|f(r+t)|+(r-t)^{m+\frac{1}{2}}|f(r-t)|\}. 
\]
It follows from (\ref{s}) that 
\begin{align*}
|J_{10}(r,t)|
&\le Cr^{m+1}\{(r+t)^{m-1}|f(r+t)|+(r-t)^{m-1}|f(r-t)|\}\le Cr^{m+1}
\end{align*}
for $0<r\le k$ and 
\begin{align*}
|J_{10}(r,t)|
&\le Cr^{m-\frac{1}{2}}\\
\end{align*}
for $r\ge k$. 
Finally we consider $D_\alpha J_{10}$. 
Differentiating $J_{10}$, we have 
\begin{align*}
|D J_{10}(r,t)|
&\le \{2m(r+t)^{2m}|f(r+t)|+(r+t)^{2m+1}|f'(r+t)|\}|K_m(r+t,r,t)|\\
&\quad+(r+t)^{2m+1}|f(r+t)|\{|(D_\lambda K_m)(r+t,r,t)|+|D K_m(r+t,r,t)|\}\\
&\quad+\{2m(r-t)^{2m}|f(r-t)|+(r-t)^{2m+1}|f'(r-t)|\}|K_m(r-t,r,t)|\\
&\quad+(r-t)^{2m+1}|f(r-t)|\{|(D_\lambda K_m)(r-t,r,t)|+|D K_m(r-t,r,t)|\}. 
\end{align*}
Using (\ref{4.30}) for $D_\lambda K_m(r+t,r,t)$, $D K_m(r+t,r,t)$ and (\ref{4.18}) with $\eta=0$ for $D_\lambda K_m(r-t,r,t)$, $D K_m(r-t,r,t)$, we have 
\begin{align*}
|D_\alpha J_{10}(r,t)|
&\le Cr^{m-\frac{1}{2}}\{(r+t)^{m-\frac{1}{2}}|f(r+t)|+(r+t)^{m+\frac{1}{2}}|f'(r+t)|\\
&\qquad\qquad\quad+(r-t)^{m-\frac{1}{2}}|f(r-t)|+(r-t)^{m+\frac{1}{2}}|f'(r-t)|\}
\end{align*}
Treating this by the same way as $J_{10}$, we have the desired estimates. 
\par
For (\ref{decay3_g}), we use (\ref{theta12}) for $|\beta|=0$ and 
and the derivative of (\ref{theta12}) for $|\beta|=1$. 
By $r>t$, $J_2\equiv 0$. 
So, we shall estimate $J_1$ and $D J_1$. 
Using $(\ref{4.4})_m$ with $\mu=0$, we have 
\[
|J_1(r,t)|
\le Cr^{m-\frac{1}{2}}\int_{r-t}^{r+t}\lambda^{m+1}|g(\lambda)|d\lambda. 
\]
It follows from (\ref{s}) that 
\[
|J_1(r,t)|
\le Cr^{m+1}\int_{r-t}^{r+t}\lambda^{m-\frac{1}{2}}|g(\lambda)|d\lambda
\le Cr^{m+1}
\]
for $r\le k$, while we simply have for $r\ge k$ 
\[
|J_1(r,t)|
\le Cr^{m-\frac{1}{2}}.
\]
\par
Differentiating $J_1$, we have 
\begin{align*}
|D_\alpha J_1(r,t)|
&\le\int_{r-t}^{r+t}\lambda^{2m+1}|g(\lambda)||(D_\alpha K_m)(\lambda,r,t)|d\lambda\\
&\quad+(r+t)^{2m+1}|g(r+t)||K_m(r+t,r,t)|+(r-t)^{2m+1}|g(r-t)||K_m(r-t,r,t)|. 
\end{align*}
Hence using (\ref{4.30}) for the second line and (\ref{4.18}) with $\mu=0$ for the others, we have 
\begin{align*}
|D_\alpha J_1(r,t)|
&\le Cr^{m-\frac{1}{2}}\left\{\int_{r-t}^{r+t}\lambda^{m-\frac{1}{2}}|g(\lambda)|d\lambda+(r+t)^{m-\frac{1}{2}}|g(r+t)|+(r-t)^{m-\frac{1}{2}}|g(r-t)|\right\}.
\end{align*}
Then, we can treat this by the same way as $J_1$ and get the desired estimate. 
This completes the proof of Proposition \ref{prop:free}. \hfill $\Box$

\section{A priori estimates}
\label{sec:apriori}

\par  
In this section, we consider the nonlinear problems (\ref{I_rad}). 
Here and hereafter, we assume $p\le q$ without loss of generality by the symmetricity of the systems. 
We will construct a solution to (\ref{I_rad}) on $\Omega_T:=(0,\infty)\times[0,T]$ in a function space $X$ defined by 
\begin{align*}
X&:=\left\{\right.(u,v)\in(C((0,\infty)\times[0,T]))^2\ :\ (\de_ru,\de_rv)\in(C((0,\infty)\times[0,T]))^2\\
&\qquad\qquad\qquad \supp u\cup\supp v\subset\{r\le t+k\},\quad \|u\|_1+\|v\|_2<\infty\left.\right\},
\end{align*}
where 
\begin{align}
\|u\|_1:=\sum_{i=0}^1\sup_{(r,t)\in \Omega_T}w_i(r,t)|D_r^iu(r,t)|,\label{norm1}\\
\|v\|_2:=\sum_{i=0}^1\sup_{(r,t)\in \Omega_T}z_i(r,t)|D_r^iz(r,t)|. \label{norm2}
\end{align}
The weight functions $w_i$, $z_i$ are defined by 
\begin{align}
w_i(r,t)&:=\dy\left(\frac{r}{k}\right)^{m-1+i}\left(\frac{r+2k}{k}\right)^{1-i}\tau_+^{\frac{1}{2}}W(r,t),\label{w}\\
z_i(r,t)&:=\dy\left(\frac{r}{k}\right)^{m-1+i}\left(\frac{r+2k}{k}\right)^{1-i}\tau_+^{\frac{1}{2}}Z(r,t)\label{z}\\
\end{align}
for $i=0,1$, where 
\begin{align}
W(r,t)&:=\tau_-^{(m+\frac{1}{2})p-(m+\frac{3}{2})},\label{W}\\
Z(r,t)&:=\left\{\begin{array}{ll}
\tau_-^{\frac{1}{p}}\left(\log 3\tau_-\right)^\nu &\mbox{if}\quad F=0,\ p\neq q\\
\tau_-^\mu & \mbox{otherwise},
\end{array}\right.\label{Z}
\end{align}
with 
\begin{align}
\mu&:=\frac{1}{p}-(q-1)F(p,q;n),\label{mu}\\
\nu&:=\frac{q(p-1)}{p(pq-1)}. \label{nu}
\end{align}
\par
Let us consider the sequences $\{U_j,V_j\}$ defined by 
\begin{equation}
\label{seq}
\begin{array}{l}
U_1=0,\qquad U_{j+1}=L(|v_0+V_j|^p)\quad (j\ge 1),\\
V_1=0,\qquad V_{j+1}=L(|u_0+U_j|^q)\quad (j\ge 1),
\end{array}
\end{equation}
where $u_0$, $v_0$ and $L$ are defined in (\ref{u0}), (\ref{v0}) and (\ref{L}), respectively. 
We will show the convergence of these sequences in $X$. 
If we prove the existence of the limit functions $U$ and $V$, 
we can show the existence of solution $u$ and $v$ to (\ref{I_rad}), 
since $u=\ep u^0+U$, $v=\ep v^0+V$. 
\par
First, recall the representation formula of $L(G)$ in (\ref{L}) and the expressions of $\Theta$  (\ref{theta12}), (\ref{theta34}), (\ref{theta56}). 
Then, we have 
\begin{align}
r^{2m}L(G)(r,t)
&=\int_0^t\left\{\w{J}_1(G(\cdot,\tau))(r,t-\tau)+\w{J}_2(G(\cdot,\tau))(r,t-\tau)\right\}d\tau, \label{L12}\\
2r^{2m}L(G)(r,t)
&=\int_0^t\left\{\w{J}_3(G(\cdot,\tau))(r,t-\tau)+\w{J}_4(G(\cdot,\tau))(r,t-\tau)\right\}d\tau, \label{L34}\\
\de_r(2r^{2m}L(G)(r,t))
&=\int_0^t\left\{\w{J}_5(G(\cdot,\tau))(r,t-\tau)+\w{J}_6(G(\cdot,\tau))(r,t-\tau)\right\}d\tau, \label{L56}
\end{align}
where $\w{J}_i$ ($i=1,2,\cdots,6$) is defined by replacing $g$ by $G(\cdot,\tau)$ in 
(\ref{J1}), (\ref{J2}), (\ref{J3}), (\ref{J4}), (\ref{J5}), (\ref{J6}).  
According to the basic attitude and the choice of the estimates of $K$, $\w{K}$ in the proof of Lemma 6.3, Lemma 6.4, Lemma 6.5 of \cite{KK_even}, we have 
\begin{align}
|L(G)(r,t)|
&\le\dy\frac{C}{r^{m+1/2-\eta}}\int_0^td\tau\int_{|\lambda_-|}^{\lambda_+}\frac{\lambda^{m+1-\eta}|G(\lambda,\tau)|}{\sqrt{\tau+\lambda-t+r}}d\lambda\nonumber\\
&\quad+\dy\frac{C}{r^{m+1/2-\eta}}\int_0^{(t-r)_+}(t-r-\tau)^{-m-\eta}d\tau
\int_0^{\lambda_-}\frac{\lambda^{2m+1}|G(\lambda,\tau)|}{\sqrt{t-r-\tau-\lambda}}d\lambda, 
\label{rep1}
\end{align}
by (\ref{L12}) and 
\begin{align}
\label{rep2}
|\de_r^iL(G)(r,t)|
&\le\dy\frac{C}{r^{m-1/2-\eta+i}}\int_0^td\tau\int_{|\lambda_-|}^{\lambda_+}\frac{|\de_\lambda(\lambda^{2m}G(\lambda,\tau))|\lambda^{-m+1-\eta}}{\sqrt{\tau+\lambda-t+r}}d\lambda\nonumber\\
&\quad+\dy\frac{C}{r^{m-1/2-\eta+i}}\int_0^{(t-r)_+}(t-r-\tau)^{m+1/2-\eta}|G(\lambda_-/2,\tau)|d\tau\nonumber\\
&\quad+\dy\frac{C}{r^{m-1/2-\eta+i}}\int_0^{(t-r)_+}(t-r-\tau)^{-m+1-\eta}d\tau\int_0^{\lambda_-/2}\frac{\lambda^{2m}|G(\lambda,\tau)|}{(t-r-\tau-\lambda)^{3/2}}d\lambda\nonumber\\
&\quad+\dy\frac{C}{r^{m-1/2-\eta+i}}\int_0^{(t-r)_+}(t-r-\tau)^{-m+1-\eta}d\tau\int_{\lambda_-/2}^{\lambda_-}\frac{\de_\lambda(\lambda^{2m}G(\lambda,\tau))|}{\sqrt{\lambda_--\lambda}}d\lambda\nonumber\\
&\quad+\dy\frac{C}{r^{m-1/2-\eta+i}}\dy\int_{t-r}^t(\tau+r-t)^{m+1/2-\eta}|G(\tau+r-t,\tau)|d\tau
\end{align}
by (\ref{L34}) and (\ref{L56}), where $\eta=0$ or $\frac{1}{2}$, $\lambda_{\pm}=t-\tau\pm r$ and $(a)_+=\max\{a,0\}$. 
\par 
To show that the sequence $(\{U_j,V_j\})$ in (\ref{seq}) converges in $X$, it is necessary to consider (\ref{rep1}) and (\ref{rep2}) with $G=G_1,\ G_2,\ G_3,\ G_4$ below.
\begin{align*}
G_1&:=|v_0+V_j|^p,\\
G_2&:=|u_0+U_j|^q,\\
G_3&:=|v_0+V_j|^p-|v_0+V_{j-1}|^p,\\
G_4&:=|u_0+U_j|^q-|u_0+U_{j-1}|^q. 
\end{align*}
First, it holds that 
\begin{align*}
|G_1|&\le C(|v_0|^p+|V_j|^p),\\
|G_2|&\le C(|u_0|^q+|U_j|^q),\\
|G_3|&\le C|V_j-V_{j-1}|(|v_0|^{p-1}+\w{V}_j^{p-1}),\\
|G_4|&\le C|U_j-U_{j-1}|(|u_0|^{q-1}+\w{U}_j^{q-1}),
\end{align*}
and 
\begin{align*}
|\de_\lambda G_1|&\le C(|v_0|^{p-1}+|V_j|^{p-1})(|\de_\lambda v_0|+|\de_\lambda V_j|),\\
|\de_\lambda G_2|&\le C(|u_0|^{q-1}+|U_j|^{q-1})(|\de_\lambda u_0|+|\de_\lambda U_j|),\\
|\de_\lambda G_3|&\le 
C|\de_\lambda V_j-\de_\lambda V_{j-1}|(|v_0|^{p-1}+|V_j|^{p-1})
+C(|\de_\lambda v_0|+|\de_\lambda V_{j-1}|)|V_j-V_{j-1}|^{p-1},\\
|\de_\lambda G_4|&\le C|\de_\lambda U_j-\de_\lambda U_{j-1}|(|u_0|^{q-1}+|U_j|^{q-1})
+C(|\de_\lambda u_0|+|\de_\lambda U_{j-1}|)|U_j-U_{j-1}|^{q-1}, 
\end{align*}
where $\w{V}_j=\max\{|V_j|,\ |V_{j-1}|\}$, $\w{U}_j=\max\{|U_j|,\ |U_{j-1}|\}$. 
Here the second terms of $\de_\lambda G_3$ and $\de_\lambda G_4$ come from H\"{o}lder continuity,  since $1<p,q<2$. 
\par 
To close the iteration, we introduce the auxiliary norm $\-\cdot\-$ as
\begin{align}
\-u\-_1:=\sup_{(r,t)\in\Omega_T}w_1(r,t)|u(r,t)|,\label{norm1'}\\
\-v\-_2:=\sup_{(r,t)\in\Omega_T}z_1(r,t)|v(r,t)|. \label{norm2'}
\end{align}
according to \cite{KK_even}. 
Hence, Using Proposition \ref{prop:free}, (\ref{norm1}), (\ref{norm2}), (\ref{norm1'}) and (\ref{norm2'}), we get the followings by (\ref{rep1}). 
\begin{align}
|U_{j+1}|
&\le\frac{C}{r^{m+\frac{1}{2}-\eta}}
\left[\ep^p\sum_{k=1}^2I_k(\phi_0^{-p},p)+\|V_j\|_2^p\sum_{k=1}^2I_k(Z^{-p},p)\right]\label{U_b}\\
|V_{j+1}|
&\le\frac{C}{r^{m+\frac{1}{2}-\eta}}
\left[\ep^q\sum_{k=1}^2I_k(\phi_0^{-q},q)+\|U_j\|_1^q\sum_{k=1}^2I_k(W^{-q},q)\right]\label{V_b}\\
|U_{j+1}-U_j|
&\le\frac{C}{r^{m+\frac{1}{2}-\eta}}\|V_j-V_{j-1}\|_2
\left[\ep^{p-1}\sum_{k=1}^2I_k(\phi_0^{1-p}Z^{-1},p)
+\|\w{V}_j\|_2^{p-1}\sum_{k=1}^2I_k(Z^{-p},p)\right]\label{U_d1}\\
|V_{j+1}-V_j|
&\le\frac{C}{r^{m+\frac{1}{2}-\eta}}\|U_j-U_{j-1}\|_1
\left[\ep^{q-1}\sum_{k=1}^2I_k(\phi_0^{1-q}W^{-1},q)
+\|\w{U}_j\|_1^{q-1}\sum_{k=1}^2I_k(W^{-q},q)\right]\label{V_d1}\\
|U_{j+1}-U_j|
&\le\frac{C}{r^{m+\frac{1}{2}-\eta}}\-V_j-V_{j-1}\-_2
\left[\ep^{p-1}\sum_{k=1}^2I_k(\phi_0^{1-p}Z^{-1},p)
+\|\w{V}_j\|_2^{p-1}\sum_{k=1}^2I_k(Z^{-p},p)\right]\label{U_d2}\\
|V_{j+1}-V_j|
&\le\frac{C}{r^{m+\frac{1}{2}-\eta}}\-U_j-U_{j-1}\-_1
\left[\ep^{q-1}\sum_{k=1}^2I_k(\phi_0^{1-q}W^{-1},q)
+\|\w{U}_j\|_1^{q-1}\sum_{k=1}^2I_k(W^{-q},q)\right], \label{V_d2}
\end{align}
where
\begin{align}
I_1(\Phi,p)(r,t)
&:=
\int_{-k}^{t-r}d\beta
\int_{|t-r|}^{t+r}
\frac{
\left(\frac{\alpha-\beta}{k}\right)^{m+1-mp-\eta}\left(\frac{\alpha+2k}{k}\right)^{-p/2}\Phi\left(\frac{\alpha-\beta}{2},\frac{\alpha+\beta}{2}\right)}
{\sqrt{\alpha-t+r}}d\alpha\label{I1}\\
I_2(\Phi,p)(r,t)
&:=
\int_{0}^{t-r}
\frac{\left(\frac{\alpha+2k}{k}\right)^{-p/2}}{\sqrt{t-r-\alpha}}d\alpha\\
&\qquad\times\int_{-k}^{\alpha}
\left(\frac{\alpha-\beta}{k}\right)^{2m-mp}\left(t-r-\frac{\alpha+\beta}{2}\right)^{-m+1-\eta}
\Phi\left(\frac{\alpha-\beta}{2},\frac{\alpha+\beta}{2}\right)d\beta, \label{I2}
\end{align}
provided $\supp U_j\cup\supp V_j\subset\{r\le t+k\}$ ($j\ge 0$). 
Similarly, it follows from (\ref{rep2}) that 
\begin{align}
|D_{r,t}^iU_{j+1}|
&\le\frac{C}{r^{m-\frac{1}{2}-\eta+i}}\\
&\qquad\times
\left[
\ep^p\sum_{k=1}^4I_k(\phi_0^{-p},p)+\ep^p\sum_{k=1}^2I_k(\phi_0^{1-p}\phi_1^{-1},p)
+\ep^{p-1}\|V_j\|_2\sum_{k=1}^2I_k(\phi_0^{1-p}Z^{-1},p)\right.\nonumber\\
&\qquad\qquad\qquad\left.
+\ep\|V_j\|_2^{p-1}\sum_{k=1}^2I_k(\phi_0^{-1}Z^{1-p},p)
+\|V_j\|_2^{p}\sum_{k=1}^4I_k(Z^{-p},p)
\right],\label{DU_b}
\end{align}
\begin{align}
|D_{r,t}^iV_{j+1}|
&\le\frac{C}{r^{m-\frac{1}{2}-\eta+i}}\\
&\qquad\times\left[
\ep^q\sum_{k=1}^4I_k(\phi_0^{-q},q)+\ep^q\sum_{k=1}^2I_k(\phi_0^{1-q}\phi_1^{-1},q)
+\ep^{q-1}\|U_j\|_1\sum_{k=1}^2I_k(\phi_0^{1-q}W^{-1},q)\right.\nonumber\\
&\qquad\qquad\qquad\left.
+\ep\|U_j\|_1^{q-1}\sum_{k=1}^2I_k(\phi_0^{-1}W^{1-q},q)
+\|U_j\|_1^{q}\sum_{k=1}^4I_k(W^{-q},q)
\right].\label{DV_b}
\end{align}
\begin{align}
&\hspace{-6pt}|D_{r,t}^iU_{j+1}-D_{r,t}^iU_j|\nonumber\\
&\le \frac{C}{r^{m-\frac{1}{2}-\eta+i}}\|V_j-V_{j-1}\|_2
\left[
\ep^{p-1}\sum_{k=1}^4I_k(\phi_0^{1-p}Z^{-1},p)+\|V_j\|_2^{p-1}\sum_{k=1}^4I_k(Z^{-p},p)
\right]\nonumber\\
&\quad
+\frac{C}{r^{m-\frac{1}{2}-\eta+i}}\-V_j-V_{j-1}\-_2^{p-1}\left[
\ep\sum_{k=1}^2I_k(\phi_1^{-1}Z^{1-p},p)+\|V_{j-1}\|_2\sum_{k=1}^2I_k(Z^{-p},p)
\right],\label{DU_d}\\
&\hspace{-6pt}|D_{r,t}^iV_{j+1}-D_{r,t}^iV_j|\nonumber\\
&\le \frac{C}{r^{m-\frac{1}{2}-\eta+i}}
\|U_j-U_{j-1}\|_1
\left[
\ep^{q-1}\sum_{k=1}^4I_k(\phi_0^{1-q}W^{-1},q)+\|U_j\|_1^{q-1}\sum_{k=1}^4I_k(W^{-q},q)
\right]\nonumber\\
&\quad
+\frac{C}{r^{m-\frac{1}{2}-\eta+i}}\-U_j-U_{j-1}\-_1^{q-1}\left[
\ep\sum_{k=1}^2I_k(\phi_1^{-1}W^{1-q},q)+\|U_{j-1}\|_1\sum_{k=1}^2I_k(W^{-q},q)
\right],\label{DV_d}
\end{align}
where
\begin{align}
I_3(\Phi,p)(r,t)
&:=\int_{\frac{t-r-2k}{3}}^{t-r}\left(\frac{\lambda_-}{k}\right)^{m+\frac{1}{2}-(m-1)p-\eta}\left(\frac{\lambda_-+2k}{k}\right)^{-p}\left(\frac{\tau+\frac{\lambda_-}{2}+2k}{k}\right)^{-p/2}\Phi\left(\frac{\lambda_-}{2}, \tau\right)d\tau\label{I3}\\
I_4(\Phi,p)(r,t)
&:=\int_{t-r}^t\left(\frac{-\lambda_-}{k}\right)^{m+\frac{1}{2}-(m-1)p-\eta}\left(\frac{-\lambda_-+2k}{k}\right)^{-p}\left(\frac{\tau-\lambda_-+2k}{k}\right)^{-p/2}\Phi(-\lambda_-,\tau)d\tau, \label{I4}
\end{align}
provided $\supp U_j\cup\supp V_j\subset\{r\le t+k\}$ ($j\ge 0$). 
We have the following estimates for $I_k$ ($k=1,2,3,4$). 
\begin{prop}
\label{prop:12}
Let $m\ge 2$, $\frac{2m+3}{2m+1}<p,q<\frac{2m+5}{2m+1}$ and $F(p,q;n)\ge 0$. Then, we have 
\begin{align*}
|I_k(\phi_0^{-p},p)|&\le C\tau_-^{-\eta}W(r,t)^{-1},&|I_k(\phi_0^{-q},q)|&\le C\tau_-^{-\eta}Z(r,t)^{-1},\\
|I_k(\phi_0^{1-p}\phi_1^{-1},p)|&\le C\tau_-^{-\eta}W(r,t)^{-1},&|I_k(\phi_0^{1-q}\phi_1^{-1},q)|&\le C\tau_-^{-\eta}Z(r,t)^{-1},\\
|I_k(\phi_0^{1-p}Z^{-1},p)|&\le C\tau_-^{-\eta}W(r,t)^{-1},&|I_k(\phi_0^{1-q}W^{-1},q)|&\le C\tau_-^{-\eta}Z(r,t)^{-1},\\
|I_k(\phi_0^{-1}Z^{1-p},p)|&\le C\tau_-^{-\eta}W(r,t)^{-1},&|I_k(\phi_0^{-1}W^{1-q},q)|&\le C\tau_-^{-\eta}Z(r,t)^{-1},\\
|I_k(\phi_1^{-1}Z^{1-p},p)|&\le C\tau_-^{-\eta}W(r,t)^{-1},&|I_k(\phi_1^{-1}W^{1-q},q)|&\le C\tau_-^{-\eta}Z(r,t)^{-1},\\
|I_k(Z^{-p},p)|&\le C\tau_-^{-\eta}W(r,t)^{-1}E_1(T),&|I_k(W^{-q},q)|&\le C\tau_-^{-\eta}Z(r,t)^{-1}E_2(T),
\end{align*}
in $\Omega_T$ for $\eta=0,\frac{1}{2}$ and $k=1,2$, where 
\begin{equation}
\label{E1}
E_1(T):=\left\{\begin{array}{cl}
\left(\log\frac{T+2k}{k}\right)^{1-p\nu}&\mbox{if}\quad F(p,q;n)=0,\ p\neq q,\\
\log\frac{T+2k}{k}&\mbox{if}\quad F(p,q;n)=0,\ p=q,\\
\left(\frac{T+2k}{k}\right)^{p(q-1)F(p,q;n)}&\mbox{if}\quad F(p,q;n)>0,\\
\end{array}\right.,
\end{equation}
\begin{equation}
\label{E2}
E_2(T):=\left\{\begin{array}{cl}
\left(\log\frac{T+2k}{k}\right)^{\nu}&\mbox{if}\quad F(p,q;n)=0,\ p\neq q,\\
\log\frac{T+2k}{k}&\mbox{if}\quad F(p,q;n)=0,\ p=q,\\
\left(\frac{T+2k}{k}\right)^{q(p-1)F(p,q;n)}&\mbox{if}\quad F(p,q;n)>0,\\
\end{array}\right..
\end{equation}
\end{prop}
\begin{prop}
\label{prop:34}
Let $m\ge 2$, $\frac{2m+3}{2m+1}<p,q<\frac{2m+5}{2m+1}$ and $F(p,q;n)\ge 0$. Then, we have 
\begin{align*}
|I_k(\phi_0^{-p},p)|&\le C\tau_-^{-\eta}W(r,t)^{-1},&|I_k(\phi_0^{-q},q)|&\le C\tau_-^{-\eta}Z(r,t)^{-1},\\
|I_k(\phi_0^{1-p}\phi_1^{-1},p)|&\le C\tau_-^{-\eta}W(r,t)^{-1},&|I_k(\phi_0^{1-q}\phi_1^{-1},q)|&\le C\tau_-^{-\eta}Z(r,t)^{-1},\\
|I_k(Z^{-p},p)|&\le C\tau_-^{-\eta}W(r,t)^{-1}E_1(T),&|I_kS(W^{-q},q)|&\le C\tau_-^{-\eta}Z(r,t)^{-1}E_2(T),
\end{align*}
in $\Omega_T$ for $\eta=0,\frac{1}{2}$ and $k=3,4$, where $E_1(T)$ and $E_2(T)$ are defined in (\ref{E1}) and (\ref{E2}), respectively. 
\end{prop}
\par 
Proposition \ref{prop:12} and Proposition \ref{prop:34} will be proved in $\S 6$. 
Here, we shall use two Propositions and proceed to a priori estimates. 
\par
In (\ref{U_b})-(\ref{V_d2}) together with Proposition \ref{prop:12}, we can find $r^{-m-\frac{1}{2}+\eta}\tau_-^{-\eta}$ in all the inequalities. 
Similarly, (\ref{DU_b})-(\ref{DV_d}) with Proposition \ref{prop:12}, Proposition \ref{prop:34} has $r^{-m+\frac{1}{2}+\eta-i}$ in common. 
For them, we can treat as follows. 
\begin{lm}
\label{r}
Let $(r,t)\in\Omega_T$ with $t-r\ge -k$. 
Assume that $I$, $J_i$ ($i=0,1$) satisfy 
\begin{align}
I&\le r^{-m-\frac{1}{2}+\eta}\tau_-^{-\eta} &\mbox{for}&\quad \eta=0,\ \frac{1}{2}\label{I}\\ 
J_i&\le r^{-m+\frac{1}{2}+\eta-i}\tau_-^{-\eta} &\mbox{for}&\quad \eta=0,\ \frac{1}{2},\quad i=0,1\label{J}
\end{align}
for $\eta=0,\ \frac{1}{2}$. 
Then, we have 
\begin{align}
I&\le C\left(\frac{r}{k}\right)^{-m+1}\left(\frac{r+2k}{k}\right)^{-1}\tau_+^{-\frac{1}{2}}
&\mbox{for\quad}& r\ge k,\label{I'}\\
I&\le C\left(\frac{r}{k}\right)^{-m}\tau_+^{-\frac{1}{2}}
&\mbox{for\quad}&  r>0,\label{I''}\\
J_i&\le C\left(\frac{r}{k}\right)^{-m+1-i}\left(\frac{r+2k}{k}\right)^{-1+i}\tau_+^{-\frac{1}{2}}
&\mbox{for\quad}& ``i=0\ and r\le k''\ \mbox{or}\ i=1.\label{J'}
\end{align}
\end{lm}
{\bf Proof.} 
Taking $\eta$ suitably as in the table below, one can easily prove them. \\
\begin{tabular}{|c|c|c|c|c|}
\hline
Area of $(r,t)$ & $r\ge k,\ t\le 2r$ & $r\ge k,\ t\ge 2r$ & $r\le k,\ t\le 2r$ & $r\le k,\ t\ge 2r$\\
\hline
\lower1.5ex\hbox{features} & \lower1.5ex\hbox{$\frac{9r}{k}\ge\frac{3(r+2k)}{k}\ge\tau_+$} & $\frac{3r}{k}\ge\frac{r+2k}{k}$ & $1\le \tau_-\le\tau_+$ & $2\le\frac{r+2k}{k}\le 3$\\
&&$3\tau_-\ge\tau_+$ & $\le \frac{3(r+2k)}{k}\le 9$ &  $3\tau_-\ge\tau_+$\\
\hline\hline
For (\ref{I'})&\lower1.5ex\hbox{ (\ref{I}), $\eta=0$} &\lower1.5ex\hbox{(\ref{I}), $\eta=\frac{1}{2}$} &-&-\\
\cline{1-1}\cline{4-5}
For (\ref{I''})&&& (\ref{I}), $\eta=\frac{1}{2}$ & (\ref{I}), $\eta=\frac{1}{2}$\\
\hline 
For (\ref{J'}), $j=0$ &-&-& \lower1.5ex\hbox{(\ref{J}), $\eta=\frac{1}{2}$} & \lower1.5ex\hbox{(\ref{J}), $\eta=\frac{1}{2}$}\\
\cline{1-3}
For (\ref{J'}), $j=1$ & (\ref{J}), $\eta=0$ & (\ref{J}), $\eta=\frac{1}{2}$ &&\\
\hline
\end{tabular}\\
\hfill $\Box$\\
\par
Now we are in a position to state a priori estimates which will give us the local existence of solution and the lower bound of the lifespan. 
For the norms (\ref{norm1}) and (\ref{norm2}), 
we shall apply 
When $r\ge k$, we use (\ref{U_b}), (\ref{V_b}), (\ref{U_d1}) and (\ref{V_d1}) and apply Proposition \ref{prop:12}, (\ref{I'}) and the definition of weight functions (\ref{w}), (\ref{z}) to them. 
For ``$i=0$ and $r\le k$'' or $i=1$ in (\ref{DU_b})(\ref{DV_b})(\ref{DU_d})(\ref{DV_d}), 
we shall apply Proposition \ref{prop:12}, Proposition \ref{prop:34} and (\ref{J'}) to them. 
As to the auxiliary norms (\ref{norm1'}), (\ref{norm2'}), we use (\ref{U_d2})(\ref{V_d2}) with Proposition \ref{prop:12}, 
(\ref{I''}). 
Then, we have the following. 
\begin{prop}
\label{prop:apriori}
Let $(U_j,V_j)\in X$ be the one in (\ref{seq}) for $j\ge 1$. 
Assume that $\frac{2m+3}{2m+1}<p\le q<\frac{2m+5}{2m+1}$ and $F(p,q;n)\ge 0$.  
Then, the following inequalities hold. 
\begin{align}
\|U_{j+1}\|_1
&\le A\ep^p+C\ep^{p-1}\|V_j\|_2+C\ep\|V_j\|_2^{p-1}+C\|V_j\|_2^pE_1(T),\label{AU_b}\\
\|V_{j+1}\|_2
&\le A\ep^q+C\ep^{q-1}\|U_j\|_1+C\ep\|U_j\|_1^{q-1}+C\|U_j\|_1^qE_2(T),\label{AV_b}\\
\|U_{j+1}-U_j\|_1
&\le C\|V_j-V_{j-1}\|_2(\ep^{p-1}+\|\w{V}_j\|_2^{p-1}E_1(T))\nonumber\\
&\quad +C\-V_j-V_{j-1}\-_2^{p-1}(\ep+\|V_{j-1}\|_2E_1(T)),\label{AU_d1}\\
\|V_{j+1}-V_j\|_2
&\le C\|U_j-U_{j-1}\|_1(\ep^{q-1}+\|\w{U}_j\|_1^{q-1}E_2(T))\nonumber\\
&\quad +C\-U_j-U_{j-1}\-_1^{q-1}(\ep+\|U_{j-1}\|_1E_2(T)),\label{AV_d1}\\
\-U_{j+1}-U_j\-_1
&\le C\-V_j-V_{j-1}\-_2(\ep^{p-1}+\|\w{V}_j\|_2^{p-1}E_1(T)),\label{AU_d2}\\
\-V_{j+1}-V_j\-_2
&\le C\-U_j-U_{j-1}\-_1(\ep^{q-1}+\|\w{U}_j\|_1^{q-1}E_2(T)),\label{AV_d2}
\end{align}
where $A\ge 1$ and $C\ge 1$ are constants independent of $\ep$. 
\end{prop}

\section{Construction of a solution and Lifespan}
\label{sec:lifespan}
\par
Theorem \ref{main} can be proved by the iteration methods using Proposition \ref{prop:apriori}. 
For this, we shall define a closed subspace  $Y$ of $X$ as follows. 
\[
Y:=\{(U,V)\in X\ :\ \|U\|_1\le 2A\ep^p,\ \|V\|_2\le 2A\ep^q\}, 
\]
where $A$ is the one in Proposition \ref{prop:apriori}. 
We will find the conditions which gurantee that $\{U_j,V_j\}$ is a convergent sequence in $Y$. 
\par
First, we note that 
\[
(U_1,V_1)=(0,0)\in Y. 
\]
Assume that $\{(U_j,V_j)\}\subset Y$. 
Then, we have 
\begin{gather}
(U_{j+1},V_{j+1}),\ (\de_rU_{j+1},\de_rV_{j+1})\in (C((0,\infty)\times[0,T]))^2,\\
\supp U_{j+1}\cup \supp V_{j+1}\subset\{r\le t+k\},
\end{gather}
by the definition of the sequences (\ref{seq}), Proposition \ref{prop:free} and the representation formulas 
(\ref{L34}), (\ref{L56}) with (\ref{J3}), (\ref{J4}), (\ref{J5}), (\ref{J6}). 
Moreover, we have by (\ref{AU_b}), (\ref{AV_b}), 
\begin{align*}
\|U_{j+1}\|_1
&\le A\ep^p+C\ep^{p-1}2A\ep^q+C\ep(2A\ep^q)^{p-1}+C(2A\ep^q)^pE_1(T),\\
\|V_{j+1}\|_2
&\le A\ep^q+C\ep^{q-1}2A\ep^p+C\ep(2A\ep^p)^{q-1}+C(2A\ep^p)^qE_2(T). 
\end{align*}
Thus, if we assume that 
\begin{align}
C\ep^{p-1}2A\ep^q+C\ep(2A\ep^q)^{p-1}+C(2A\ep^q)^pE_1(T)&\le A\ep^p,\label{condi:U1}\\
C\ep^{q-1}2A\ep^p+C\ep(2A\ep^p)^{q-1}+C(2A\ep^p)^qE_2(T)&\le A\ep^q,\label{condi:V1}
\end{align}
then, $(U_{j+1},V_{j+1})\in Y$. 
By the induction, 
\begin{equation}
\label{bounded}
(U_j,V_j)\in Y\quad\mbox{for any}\ j\ge 1.
\end{equation} 
Hereafter, we assume that (\ref{condi:U1}) and (\ref{condi:V1}) are fulfilled. 
\par
Next, we consider the differences of $\{(U_m,V_m)\}$. 
As to the norm $\-\cdot\-$, 
 (\ref{AU_d2}), (\ref{AV_d2}) and (\ref{bounded}) give us  
\begin{align*}
\-U_{j+1}-U_j\-_1
&\le C\-V_j-V_{j-1}\-_2(\ep^{p-1}+(2A\ep^q)^{p-1}E_1(T)),\\
\-V_{j+1}-V_j\-_2
&\le C\-U_j-U_{j-1}\-_1(\ep^{q-1}+(2A\ep^p)^{q-1}E_2(T)). 
\end{align*}
Thus, if we assume that 
\begin{equation}
C^2\{\ep^{p-1}+C(2A\ep^q)^{p-1}E_1(T))\}
\{\ep^{q-1}+C(2A\ep^p)^{q-1}E_2(T))\}
\le\frac{1}{4}, 
\label{condi:UV1}\\
\end{equation}
then, we have 
\[
\-U_{j+1}-U_j\-_1
\le\frac{1}{4}\-U_{j-1}-U_{j-2}\-_1,\qquad
\-V_{j+1}-V_j\-_2
\le\frac{1}{4}\-V_{j-1}-V_{j-2}\-_2
\]
which means that 
\begin{equation}
\-U_{j+1}-U_j\-_1\le\frac{B}{2^j},\qquad
\-V_{j+1}-V_j\-_2\le\frac{B}{2^j},
\label{diff2}
\end{equation}
where we set $B=\max\{1,\ 4\-U_3-U_2\-_1,\ 2\-U_2-U_1\-_1,\ 4\-V_3-V_2\-_2,\ 2\-V_2-V_1\-_2\}\ge 1$. 
Hereafter, we assume (\ref{diff2}). 
\par
As to the norm $\|\cdot\|$, 
by (\ref{AU_d1}), (\ref{AU_d2}), (\ref{bounded}) and (\ref{diff2}), we have 
\begin{align*}
\|U_{j+1}-U_j\|_1
&\le C\|V_j-V_{j-1}\|_2(\ep^{p-1}+(2A\ep^q)^{p-1}E_1(T))+C\left(\frac{B}{2^{j-1}}\right)^{p-1}(\ep+2A\ep^qE_1(T)),\\
\|V_{j+1}-V_j\|_2
&\le C\|U_j-U_{j-1}\|_1(\ep^{q-1}+(2A\ep^p)^{q-1}E_2(T))+C\left(\frac{B}{2^{j-1}}\right)^{q-1}(\ep+2A\ep^pE_2(T)). 
\end{align*}
Here we set $s=2^{1-p}<1$. 
If we assume that 
\begin{align}
C^2(\ep^{p-1}+(2A\ep^q)^{p-1}E_1(T))(\ep^{q-1}+(2A\ep^p)^{q-1}E_2(T))&\le s^2\label{condi:UV2}\\
C^2B^{q-1}(\ep+2A\ep^pE_2(T))(\ep^{p-1}+(2A\ep^q)^{p-1}E_1(T))&\le\frac{s^2}{2}\label{condi:U2},\\
CB^{p-1}(\ep+2A\ep^qE_1(T))&\le \frac{s}{2}\label{condi:U3}\\
C^2B^{p-1}(\ep+2A\ep^qE_1(T))(\ep^{q-1}+(2A\ep^p)^{q-1}E_2(T))&\le\frac{s^2}{2},\label{condi:V2}\\
CB^{q-1}(\ep+2A\ep^pE_2(T))&\le\frac{s}{2}, \label{condi:V3}
\end{align}
then, we have 
\[
\|U_{j+1}-U_j\|_1
\le s^2\|U_{j-1}-U_{j-2}\|_1+s^j,\qquad
\|V_{j+1}-V_j\|_2
\le s^2\|V_{j-1}-V_{j-2}\|_2+s^j, 
\]
which mean that 
\[
\|U_{j+1}-U_j\|_1\le B's^j+js^j,\qquad
\|V_{j+1}-V_j\|_2\le B's^j+js^j, 
\]
where we set 
$B'=\max\{s^{-2}\|U_3-U_2\|_1,\ s^{-1}\|U_2-U_1\|_1,\ s^{-2}\|V_3-V_2\|_2,\ s^{-1}\|V_2-V_1\|_2\}$. 
Since the right-hand sides in the above inequalities converge, we find that $\{U_j\}$ and $\{V_j\}$ converge with respect to the norm $\|\cdot\|$. 
Hence, we now know that $\{(U_j,V_j)\}$ is a convergent sequence in $Y$, 
provided (\ref{condi:U1}), (\ref{condi:V1}), (\ref{condi:UV1}), (\ref{condi:UV2}), (\ref{condi:U2}), (\ref{condi:U3}),  (\ref{condi:V2}) and (\ref{condi:V3}) are satisfied. 
\par
Let us define  
\begin{align*}
\ep_1
&:=\min\left\{
(6C)^{-\frac{1}{p-1}},\ (2^{q-1}3CA^{q-2})^{-\frac{1}{(p-1)(q-1)}},\ (8C^2B^{p-1}s^{-2})^{-\frac{1}{q}},\ (4CB^{q-1}s^{-1})^{-1}
\right\}
\end{align*}
and 
\begin{align}
E
&:=\min\left\{\right.
(2^q3A^{q-1}C)^{-1},\ (2^{q+1}A^{q-1}C^2s^{-2})^{-1},\nonumber\\
&\qquad\qquad
(2^{p+2}A^{p-1}C^2B^{q-1}s^{-2})^{-1},\ (2^{q+2}A^{q-1}C^2B^{p-1}s^{-2})^{-1},\ (8ACB^{q-1}s^{-1})^{-1}
\left.\right\}<1. 
\label{E}
\end{align}
Then, we see that for any $\ep\in (0,\ep_1)$, all of the conditions above are satisfied, provided 
\begin{align}
\ep^{p(q-1)}E_1(T)
&\le E,\label{L1}\\
\ep^{q(p-1)}E_2(T)
&\le E.\label{L2}
\end{align}
Note that we used the fact that $\ep<\ep^{p-1}$ for $1<p<2$ and $0<\ep<1$ in the calculations above. 
\par 
Now we can get the lower bound of the lifespan by (\ref{L1}) and (\ref{L2}). 
For example, let us see the case where $F=0,\ p\neq q$. 
(\ref{L1}) and (\ref{L2}) can be rewritten by (\ref{E1}) and (\ref{E2}) as 
\begin{equation}
\label{2condi}
\ep^{p(q-1)}\left(\log\frac{T+2k}{k}\right)^{1-p\nu}
\le E,\qquad 
\ep^{q(p-1)}\left(\log\frac{T+2k}{k}\right)^{\nu}
\le E. 
\end{equation}
By the definition of $\nu$ (\ref{nu}), one can see that these are guaranteed by 
\begin{equation}
\label{1condi}
T\le k\exp\left(E^{\frac{p(pq-1)}{q(p-1)}}\ep^{-p(pq-1)}\right)-2k, 
\end{equation}
since $0<E<1$ and $p\le q$. 
Besides, defining $\ep_2$ so that 
\[
k\exp\left(E^{\frac{p(pq-1)}{q(p-1)}}\ep_2^{-p(pq-1)}\right)-2k=
\exp\left(E^{\frac{p(pq-1)}{q(p-1)}}\ep_2^{-p(pq-1)}\right), 
\]
we see that  
\[
T\le\exp\left(E^{\frac{p(pq-1)}{q(p-1)}}\ep^{-p(pq-1)}\right)
\quad\mbox{for any}\quad \ep\in(0,\min\{\ep_1,\ \ep_2\}], 
\]
is a sufficient condition for (\ref{1condi}). 
\par
Similarly, we can consider the other cases. 
Therefore, we find that there exists a positive constant $\ep_0=\ep_0(p,q,f_1,f_2,g_1,g_2,n,k)$ such that, 
for any $\ep$ satisfying $0<\ep\le \ep_0$, 
a local-in-time solution to (\ref{I_rad}) in $Y$ exists, provided 
\[
T\le\left\{\begin{array}{cl}
\exp(\w{C}\ep^{-p(pq-1)})&\mbox{if}\quad F=0,\ p\neq q,\\
\exp(\w{C}\ep^{-p(p-1)})&\mbox{if}\quad F=0,\ p=q,\\
\w{C}\ep^{-F(p,q;n)^{-1}}&\mbox{if}\quad F>0.
\end{array}\right.
\]
This implies the lower bound of the lifespan.

\section{Proof of A priori estimates}
\label{sec:prop}
\par
In this section, we shall prove Proposition \ref{prop:12} and Proposition \ref{prop:34}. 
As to Proposition \ref{prop:12}, it is enough to consider $I_k$ ($k=1,2,3,4$) only for
\begin{equation}
\begin{array}{llll}
\phi_0^{-p}, & \phi_0^{1-p}Z^{-1}, & \phi_0^{-1}Z^{1-p}, & Z^{-p},\\
\phi_0^{-q}, & \phi_0^{1-q}W^{-1}, &\phi_0^{-1}W^{1-q}, &W^{-q},
\end{array}\label{8}
\end{equation}
since $\phi_1^{-1}=\tau_-^{-(m+\frac{3}{2})}\le\tau_-^{-(m+\frac{1}{2})}=\phi_0^{-1}$. 
The calculation in this section is essentially based on the way of Proposition 6.6 in \cite{KK_even}. 
\par 
Before the proofs, we remark the properties of $\mu$ in (\ref{mu}) and $\nu$ in (\ref{nu}). 
\begin{re}
Let $F(p,q;n)$ with $\frac{2m+3}{2m+1}<p\le q<\frac{2m+5}{2m+1}$. 
Then, we have 
\begin{gather}
0<\mu<\frac{1}{p},\qquad
0<\nu<\frac{1}{p},\label{munu>0}\\
1-p\mu=p(q-1)F(p,q;n),\qquad
1-p\nu=\frac{q-1}{pq-1},\label{1-pmunu}\\
m+\frac{3}{2}-\left(m+\frac{1}{2}\right)q+\mu=\frac{(p-q)(pq+1)}{p(pq-1)}\le 0,\label{5.13}\\
-\left\{\left(m+\frac{1}{2}\right)p-\left(m+\frac{3}{2}\right)\right\}q+m+\frac{5}{2}-\left(m+\frac{1}{2}\right)q=-\mu+q(p-1)F(p,q;n).\label{5.14}
\end{gather}
One can readily prove them by the definitions of $\mu$ and $\nu$. 
\end{re}
\par 
Hereafter, we assume that 
$m\ge 2$, $\frac{2m+3}{2m+1}<p\le q<\frac{2m+5}{2m+1}$ and $F(p,q;n)\ge 0$. 
Proposition \ref{prop:12} will proved in {\bf 6.1}, {\bf 6.2} and Proposition \ref{prop:34} will proved in {\bf 6.3}, {\bf 6.4}. 
\subsection{Estimates for $I_1$}
\underline{The first line of (\ref{8})}
\par
When $F=0$ with $p\neq q$, we have by (\ref{I1}), (\ref{Z}) and $\phi_0=\tau_-^{-m-\frac{1}{2}}$, 
\begin{align}
|I_1(\phi_0^{-p},p)|
&\le\int_{-k}^{t-r}\left(\frac{\beta+2k}{k}\right)^{-(m+\frac{1}{2})p}d\beta
\int_{|t-r|}^{t+r}\frac{\left(\frac{\alpha-\beta}{k}\right)^{m+1-mp-\eta}\left(\frac{\alpha+2k}{k}\right)^{-\frac{p}{2}}}{\sqrt{\alpha-(t-r)}}d\alpha\label{I11_u1}\\
|I_1(\phi_0^{1-p}Z^{-1},p)|
&\le\int_{-k}^{t-r}\left(\frac{\beta+2k}{k}\right)^{-(m+\frac{1}{2})(p-1)-\frac{1}{p}}\left(\log 3\frac{\beta+2k}{k}\right)^{-\nu}d\beta\nonumber\\
&\qquad\qquad\times
\int_{|t-r|}^{t+r}\frac{\left(\frac{\alpha-\beta}{k}\right)^{m+1-mp-\eta}\left(\frac{\alpha+2k}{k}\right)^{-\frac{p}{2}}}{\sqrt{\alpha-(t-r)}}d\alpha\label{I11_u2}\\
|I_1(\phi_0^{-1}Z^{1-p},p)|
&\le\int_{-k}^{t-r}\left(\frac{\beta+2k}{k}\right)^{-(m+\frac{1}{2})-\frac{p-1}{p}}\left(\log 3\frac{\beta+2k}{k}\right)^{-(p-1)\nu}d\beta\nonumber\\
&\qquad\qquad\times
\int_{|t-r|}^{t+r}\frac{\left(\frac{\alpha-\beta}{k}\right)^{m+1-mp-\eta}\left(\frac{\alpha+2k}{k}\right)^{-\frac{p}{2}}}{\sqrt{\alpha-(t-r)}}d\alpha\label{I11_u3}\\
|I_1(Z^{-p},p)|
&\le\int_{-k}^{t-r}\left(\frac{\beta+2k}{k}\right)^{-1}\left(\log 3\frac{\beta+2k}{k}\right)^{-p\nu}d\beta\nonumber\\
&\qquad\qquad\times
\int_{|t-r|}^{t+r}\frac{\left(\frac{\alpha-\beta}{k}\right)^{m+1-mp-\eta}\left(\frac{\alpha+2k}{k}\right)^{-\frac{p}{2}}}{\sqrt{\alpha-(t-r)}}d\alpha.\label{I11_u4}
\end{align}
We see that the $\alpha$-integrals in four inequqlities are common. 
By 
\[
\left(\frac{\alpha+2k}{k}\right)^{-\frac{p}{2}}\le \tau_-^{-\eta}\left(\frac{\alpha-\beta}{k}\right)^{-\frac{p}{2}+\eta}
\]
and the integration by parts, we have
\begin{align*}
(\alpha\mbox{-integral})
&\le C\tau_-^{-\eta}\left\{\sqrt{r}\left(\frac{t+r-\beta}{k}\right)^{m+1-(m+\frac{1}{2})p}
+\int_{|t-r|}^{t+r}\left(\frac{\alpha-\beta}{k}\right)^{m+\frac{1}{2}-(m+\frac{1}{2})p}d\alpha\right\}\\
&\le C\tau_-^{-\eta}
\left(\frac{t-r-\beta}{k}\right)^{m+\frac{3}{2}-(m+\frac{1}{2})p}, 
\end{align*}
since $-1<m+\frac{3}{2}-(m+\frac{1}{2})p<0$. Thus, for example, (\ref{I11_u1}) becomes 
\begin{equation}
\label{I1_beta1}
|I_1(\phi_0^{-p},p)|
\le C\tau_-^{-\eta}\int_{-k}^{t-r}\left(\frac{\beta+2k}{k}\right)^{-(m+\frac{1}{2})p}\left(\frac{t-r-\beta}{k}\right)^{m+\frac{3}{2}-(m+\frac{1}{2})p}d\beta. 
\end{equation}
Here, the following Lemma is useful. 
\begin{lm}
\label{lm5.4}
Let $k>0$, $-k<a$, $l>0$ and $0<h<1$. Then we have 
\[
\int_{-k}^a\left(\frac{\beta+2k}{k}\right)^{-l}\left(\frac{a-\beta}{k}\right)^{-h}d\beta\\
\le Ck\left(\frac{a+2k}{k}\right)^{-h}\w{E}(a),
\]
where 
\[
\w{E}(a):=\left\{\begin{array}{cl}
1&\mbox{if}\quad l>1,\\
\log\frac{a+2k}{k} &\mbox{if}\quad l=1,\\
\left(\frac{a+2k}{k}\right)^{1-l}&\mbox{if}\quad l<1.
\end{array}\right.
\]
\end{lm}
{\bf Proof.} 
Lemma 4.4 in \cite{Kodd} is similar to this lemma and we can find its proof there. 
So, we omit the proof here. \hfill $\Box$\\\\
Using Lemma \ref{lm5.4} with $a=t-r$, $h=(m+\frac{1}{2})p-(m+\frac{3}{2})\in(0,1)$, $l=(m+\frac{1}{2})p>1$ for 
(\ref{I1_beta1}), we have 
\[
|I_1(\phi_0^{-p},p)|
\le C\tau_-^{-\eta+m+\frac{3}{2}-(m+\frac{1}{2})p}=C\tau_-^{-\eta}W(r,t)^{-1}.
\]
As to (\ref{I11_u2}) and (\ref{I11_u3}), removing $\log\frac{\beta+2k}{k}$ and taking $l=(m+\frac{1}{2}(p-1)+\frac{1}{p}>1$, $l=m+\frac{1}{2}+\frac{p-1}{p}>1$ respectively in applying Lemma \ref{lm5.4}, 
we have the same conclusion as (\ref{I11_u1}). 
(\ref{I11_u4}) became
needs the following Lemma.  
\begin{lm}
\label{lm5.5}
Let $k>0$, $-k<a$, $0<l<1$ and $0<h<1$. Then we have 
\[
\int_{-k}^a\left(\frac{\beta+2k}{k}\right)^{-1}\left(\log 3\frac{\beta+2k}{k}\right)^{-l}\left(\frac{a-\beta}{k}\right)^{-h}d\beta\\
\le Ck\left(\frac{a+2k}{k}\right)^{-h}\left(\log 3\frac{a+2k}{k}\right)^{1-l}.
\]
\end{lm}
{\bf Proof.} 
This lemma can be proved by the same division as the one in the proof of Lemma \ref{lm5.4}. 
So, we omit the proof here. \hfill $\Box$\\\\
Applying Lemma \ref{lm5.5} with $a=t-r$, $h=(m+\frac{1}{2})p-(m+\frac{3}{2})\in(0,1)$, $l=p\nu\in(0,1)$, we get 
\[
|I_1(Z^{-p},p)|
\le C\tau_-^{-\eta+m+\frac{3}{2}-(m+\frac{1}{2})p+1-p\nu}
\le C\tau_-^{-\eta}W(r,t)^{-1}E_1(T). 
\]
\par
When $F=0$ with $p=q$ or $F>0$, we have 
\[
|I_1(\phi_0^{-p},p)|,\ |I_1(\phi_0^{1-p}Z^{-1},p)|,\ |I_1(\phi_0^{-1}Z^{1-p},p)|
\le C\tau_-^{-\eta}W(r,t)^{-1}, 
\]
by the similar calculation. 
$I_1(Z^{-p},p)$ has the following form. 
\begin{equation}
|I_1(Z^{-p},p)|
\le\int_{-k}^{t-r}\left(\frac{\beta+2k}{k}\right)^{-p\mu}d\beta
\int_{|t-r|}^{t+r}\frac{\left(\frac{\alpha-\beta}{k}\right)^{m+1-mp-\eta}\left(\frac{\alpha+2k}{k}\right)^{-\frac{p}{2}}}{\sqrt{\alpha-(t-r)}}d\alpha.\label{I12_u4}
\end{equation}
After the same treatment for the $\alpha$-integral, we shall apply Lemma \ref{lm5.4} with 
$a=t-r$, $h=(m+\frac{1}{2})p-(m+\frac{3}{2})\in(0,1)$ and $l=p\mu$ to the $\beta$-integral. 
Noting that $l=1$ if $F(p,q;n)=0$ and $l<1$ if $F(p,q;n)>0$ because of (\ref{1-pmunu}), we have 
\begin{align*}
|I_1(Z^{-p},p)|
&\le C\tau_-^{-\eta+m+\frac{3}{2}-(m+\frac{1}{2})p}\times\left\{\begin{array}{cl}
\log \tau_-&\mbox{if}\quad F(p,q;n)=0,\ p=q\\
\tau_-^{1-p\mu}&\mbox{if}\quad F(p,q;n)>0
\end{array}\right.\\
&\le C\tau_-^{-\eta}W(r,t)^{-1}E_1(T). 
\end{align*}
\noindent
\underline{The second line of (\ref{8})}
\par
By the similar way to the above case, we have 
\[
\label{I1_v123}
|I_1(\phi_0^{-q},q)|,\  
|I_1(\phi_0^{1-q}W^{-1},q)|,\  
|I_1(\phi_0^{-1}W^{1-q},q)|
\le C\tau_-^{-\eta+m+\frac{3}{2}-(m+\frac{1}{2})q}. 
\]
Hence, we get the desired estimates by (\ref{5.13}). 
As to $I_1(W^{-q},q)$, in applying Lemma \ref{lm5.4}, we take $l=\{(m+\frac{1}{2})p-(m+\frac{3}{2})\}q$. 
Using (\ref{5.14}) and (\ref{5.13}) in turn, we see that 
\begin{align}
l&=m+\frac{5}{2}-(m+\frac{1}{2})q+\mu-q(p-1)F(p,q;n)\nonumber\\
&=1+\frac{(p-q)(pq+1)}{p(pq-1)}-q(p-1)F(p,q;n)\label{l}
\end{align}
which yields that $l=1$ if $F(p,q;n)=0$, $p=q$ and $l<1$ if $F(p,q;n)>0$. 
Thus we get 
\[
|I_1(W^{-q},q)|\le C\tau_-^{-\eta+m+\frac{3}{2}-(m+\frac{1}{2})q}
\times\left\{\begin{array}{cl}
\log\tau_-&\mbox{if}\quad F(p,q;n)=0,\ p=q\\
\tau_-^{1-\{(m+\frac{1}{2})p-(m+\frac{3}{2})\}q}&\mbox{Otherwise}. 
\end{array}\right.
\]
Using (\ref{5.13}) if $F(p,q;n)=0,\ p=q$ and (\ref{5.14}) if $F(p,q;n)>0$, this is dominated by $C\tau_-^{-\eta}Z(r,t)^{-1}E_2(T)$. 
When $F=0,\ p\neq q$, noting that $\mu=\frac{1}{p}$ by (\ref{1-pmunu}) and using (\ref{5.14}), we have 
\begin{align*}
|I_1(W^{-q},q)|
&\le C\tau_-^{-\eta-\frac{1}{p}}\\
&=C\tau_-^{-\eta-\frac{1}{p}}(\log 3\tau_-)^{-\nu}(\log 3\tau_-)^\nu\\
&\le C\tau_-^{-\eta}Z(r,t)^{-1}E_2(T).
\end{align*}
Therefore Proposition \ref{prop:12} for $I_1$ has been proved.

\subsection{Estimates for $I_2$}
\par
Here, we assume that $t-r>0$. 
We shall divide the $\alpha$-integral in (\ref{I2}) as follows. 

\begin{equation}
\int_0^{t-r}d\alpha=\int_{0}^{(\frac{t-r}{2}-k)_+}+\int_{\frac{t-r}{2}-k}^{t-r}d\alpha. 
\label{I2'}
\end{equation}
Since it holds that 
\[
\begin{array}{l}
\dy t-r-\frac{\alpha+\beta}{2}\ge t-r-\alpha,\\
\dy t-r-\frac{\alpha+\beta}{2}\ge \frac{\alpha-\beta}{2},\\
\dy t-r-\frac{\alpha+\beta}{2}\le C\tau_-
\end{array}
\]
in the domain of the integration, we have two estimates: 
\begin{align}
&\left(\frac{\alpha-\beta}{k}\right)^{2m-mp}\left(t-r-\frac{\alpha+\beta}{2}\right)^{-m+1-\eta}\nonumber\\
&\qquad\qquad\le
\dy\left\{\begin{array}{l}
\dy C\tau_-^{\frac{1}{2}-\eta}(t-r-\alpha)^{m+\frac{3}{2}-(m+\frac{1}{2})p-\delta}
\left(\frac{\alpha-\beta}{k}\right)^{\frac{p}{2}-1+\delta},\\
\dy C\tau_-^{m+\frac{5}{2}-(m+\frac{1}{2})p}(t-r-\alpha)^{-\frac{1}{2}+\delta}
\left(\frac{\alpha-\beta}{k}\right)^{\frac{p}{2}-1-\eta-\delta},
\end{array}\right.
\label{shitaue}
\end{align}
for $0<\delta<1/2$ and $\eta=0$ or $1/2$. 
Using the first line of (\ref{shitaue}) for the first term in (\ref{I2'}) and the second line of (\ref{shitaue}) for the second term in (\ref{I2'}), we have 
\begin{align}
|I_2(\Phi,p)|
&=C\tau_-^{\frac{1}{2}-\eta}\int_{0}^{(\frac{t-r}{2}-k)_+}
\left(\frac{\alpha+2k}{k}\right)^{-p/2}(t-r-\alpha)^{m+1-(m+\frac{1}{2})p-\delta}d\alpha\nonumber\\
&\quad\qquad\times\int_{-k}^{\alpha}
\left(\frac{\alpha-\beta}{k}\right)^{\frac{p}{2}-1+\delta}
\Phi\left(\frac{\alpha-\beta}{2},\frac{\alpha+\beta}{2}\right)d\beta\nonumber\\
&\quad+C\tau_-^{m+\frac{5}{2}-(m+\frac{1}{2})p}\int_{\frac{t-r}{2}-k}^{t-r}
\left(\frac{\alpha+2k}{k}\right)^{-p/2}(t-r-\alpha)^{-1+\delta}d\alpha\nonumber\\
&\quad\qquad\times\int_{-k}^{\alpha}
\left(\frac{\alpha-\beta}{k}\right)^{\frac{p}{2}-1-\eta-\delta}
\Phi\left(\frac{\alpha-\beta}{2},\frac{\alpha+\beta}{2}\right)d\beta.\label{I2''}
\end{align}
\noindent
\underline{The first line of (\ref{8})}
\par
When $F=0,\ p\neq q$, we have by (\ref{I2''}), (\ref{Z}) and $\phi_0=\tau_-^{-m-\frac{1}{2}}$, 
\begin{align}
|I_2(\phi_0^{-p},p)|
&=C\tau_-^{\frac{1}{2}-\eta}\int_{0}^{(\frac{t-r}{2}-k)_+}
\left(\frac{\alpha+2k}{k}\right)^{-p/2}(t-r-\alpha)^{m+1-(m+\frac{1}{2})p-\delta}d\alpha\nonumber\\
&\qquad\qquad\times\int_{-k}^{\alpha}
\left(\frac{\alpha-\beta}{k}\right)^{\frac{p}{2}-1+\delta}
\left(\frac{\beta+2k}{k}\right)^{-(m+\frac{1}{2})p}d\beta\nonumber\\
&\quad+C\tau_-^{m+\frac{5}{2}-(m+\frac{1}{2})p}\int_{\frac{t-r}{2}-k}^{t-r}
\left(\frac{\alpha+2k}{k}\right)^{-p/2}(t-r-\alpha)^{-1+\delta}d\alpha\nonumber\\
&\qquad\qquad\times\int_{-k}^{\alpha}
\left(\frac{\alpha-\beta}{k}\right)^{\frac{p}{2}-1-\eta-\delta}
\left(\frac{\beta+2k}{k}\right)^{-(m+\frac{1}{2})p}d\beta,\label{I2_u1}
\end{align}
\begin{align}
|I_2(\phi_0^{1-p}Z^{-1},p)|
&=C\tau_-^{\frac{1}{2}-\eta}\int_{0}^{(\frac{t-r}{2}-k)_+}
\left(\frac{\alpha+2k}{k}\right)^{-p/2}(t-r-\alpha)^{m+1-(m+\frac{1}{2})p-\delta}d\alpha\nonumber\\
&\quad\qquad\times\int_{-k}^{\alpha}
\left(\frac{\alpha-\beta}{k}\right)^{\frac{p}{2}-1+\delta}
\left(\frac{\beta+2k}{k}\right)^{-(m+\frac{1}{2})(p-1)-\frac{1}{p}}\left(\log 3\frac{\beta+2k}{k}\right)^{-\nu}d\beta\nonumber\\
&\quad+C\tau_-^{m+\frac{5}{2}-(m+\frac{1}{2})p}\int_{\frac{t-r}{2}-k}^{t-r}
\left(\frac{\alpha+2k}{k}\right)^{-p/2}(t-r-\alpha)^{-1+\delta}d\alpha\nonumber\\
&\quad\qquad\times\int_{-k}^{\alpha}
\left(\frac{\alpha-\beta}{k}\right)^{\frac{p}{2}-1-\eta-\delta}
\left(\frac{\beta+2k}{k}\right)^{-(m+\frac{1}{2})(p-1)-\frac{1}{p}}\left(\log 3\frac{\beta+2k}{k}\right)^{-\nu}d\beta,\label{I2_u2}
\end{align}
\begin{align}
|I_2(\phi_0^{-1}Z^{1-p},p)|
&=C\tau_-^{\frac{1}{2}-\eta}\int_{0}^{(\frac{t-r}{2}-k)_+}
\left(\frac{\alpha+2k}{k}\right)^{-p/2}(t-r-\alpha)^{m+1-(m+\frac{1}{2})p-\delta}d\alpha\nonumber\\
&\qquad\times\int_{-k}^{\alpha}
\left(\frac{\alpha-\beta}{k}\right)^{\frac{p}{2}-1+\delta}
\left(\frac{\beta+2k}{k}\right)^{-(m+\frac{1}{2})-\frac{p-1}{p}}\left(\log 3\frac{\beta+2k}{k}\right)^{-(p-1)\nu}d\beta\nonumber\\
&\quad+C\tau_-^{m+\frac{5}{2}-(m+\frac{1}{2})p}
\int_{\frac{t-r}{2}-k}^{t-r}
\left(\frac{\alpha+2k}{k}\right)^{-p/2}(t-r-\alpha)^{-1+\delta}d\alpha\nonumber\\
&\qquad\times\int_{-k}^{\alpha}
\left(\frac{\alpha-\beta}{k}\right)^{\frac{p}{2}-1-\eta-\delta}
\left(\frac{\beta+2k}{k}\right)^{-(m+\frac{1}{2})-\frac{p-1}{p}}\left(\log 3\frac{\beta+2k}{k}\right)^{-(p-1)\nu}d\beta,\label{I2_u3}
\end{align}
\begin{align}
|I_2(Z^{-p},p)|
&=C\tau_-^{\frac{1}{2}-\eta}\int_{0}^{(\frac{t-r}{2}-k)_+}
\left(\frac{\alpha+2k}{k}\right)^{-p/2}(t-r-\alpha)^{m+1-(m+\frac{1}{2})p-\delta}d\alpha\nonumber\\
&\qquad\qquad\times\int_{-k}^{\alpha}
\left(\frac{\alpha-\beta}{k}\right)^{\frac{p}{2}-1+\delta}
\left(\frac{\beta+2k}{k}\right)^{-1}\left(\log 3\frac{\beta+2k}{k}\right)^{-p\nu}d\beta\nonumber\\
&\quad+C\tau_-^{m+\frac{5}{2}-(m+\frac{1}{2})p}\int_{\frac{t-r}{2}-k}^{t-r}
\left(\frac{\alpha+2k}{k}\right)^{-p/2}(t-r-\alpha)^{-1+\delta}d\alpha\nonumber\\
&\qquad\qquad\times\int_{-k}^{\alpha}
\left(\frac{\alpha-\beta}{k}\right)^{\frac{p}{2}-1-\eta-\delta}
\left(\frac{\beta+2k}{k}\right)^{-1}\left(\log 3\frac{\beta+2k}{k}\right)^{-p\nu}d\beta.\label{I2_u4}
\end{align}
for $0<\delta<\frac{1}{2}$. 
\par
Let us Remove $\log\frac{\beta+2k}{k}$ in (\ref{I2_u2}), (\ref{I2_u3}). 
Then, (\ref{I2_u1}), (\ref{I2_u2}) and (\ref{I2_u3}) can be treated by the same way. 
As to the first term, setting $a=\alpha$, $h=-\frac{p}{2}+1-\delta$, 
we have $l>1$, $h\in(0,1)$ provided $0<\delta<\frac{1}{10}$. 
Thus, applying Lemma \ref{lm5.4} with $a$, $h$ in common, and $l=(m+\frac{1}{2})p>1$, $l=(m+\frac{1}{2})(p-1)+\frac{1}{p}>1$, $l=m+\frac{1}{2}+\frac{p-1}{p}>1$ in turn, we can calculate them as follows. 
\begin{align*}
&\mbox{(the first tems of \ref{I2_u1}), (\ref{I2_u2}), (\ref{I2_u3}))}\\
&\qquad \le C\tau_-^{\frac{1}{2}-\eta}
\int_0^{(\frac{t-r}{2}-k)_+}\left(\frac{\alpha+2k}{k}\right)^{-1+\delta}(t-r-\alpha)^{m+1-(m+\frac{1}{2})p-\delta}d\alpha\\
&\qquad\le C\tau_-^{-\eta+m+\frac{3}{2}-(m+\frac{1}{2})p-\delta}
\int_0^{(\frac{t-r}{2}-k)_+}\left(\frac{\alpha+2k}{k}\right)^{-1+\delta}
d\alpha\\
&\qquad\le C\tau_-^{-\eta+m+\frac{3}{2}-(m+\frac{1}{2})p}\\
&\qquad\le C\tau_-^{-\eta}W(r,t)^{-1}.
\end{align*}
As to the second term, setting $a=\alpha$ and $h=-\frac{p}{2}+1+\eta+\delta$, 
we have $h\in(0,1)$ provided $0<\delta<\frac{1}{2m+1}$. 
Taking same $l$ as the first terms and applying Lemma \ref{lm5.4}, we have 
\begin{align*}
&\mbox{(the second term of \ref{I2_u1}), (\ref{I2_u2}), (\ref{I2_u3}))}\\
&\qquad\le C\tau_-^{m+\frac{5}{2}-(m+\frac{1}{2})p}
\int_{\frac{t-r}{2}-k}^{t-r}\left(\frac{\alpha+2k}{k}\right)^{-1-\eta-\delta+\delta'}(t-r-\alpha)^{-1+\delta}d\alpha\\
&\qquad\le C\tau_-^{-\mu+m+\frac{3}{2}-(m+\frac{1}{2})p-\delta}
\int_{\frac{t-r}{2}-k}^{t-r}(t-r-\alpha)^{-1+\delta}d\alpha\\
&\qquad\le C\tau_-^{-\mu+m+\frac{3}{2}-(m+\frac{1}{2})p}\\
&\qquad= C\tau_-^{-\eta}W(r,t)^{-1}.
\end{align*}
Therefore we get 
\[
|I_2(\phi_0^{-p},p)|,\quad |I_2(\phi_0^{1-p}Z^{-1},p)|,\quad |I_2(\phi_0^{-1}Z^{1-p},p)|\le C\tau_-^{-\eta}W(r,t)^{-1}.\\
\]
As to (\ref{I2_u4}), we do not remove $\log\frac{\beta+2k}{k}$ and apply Lemma \ref{lm5.5}. 
Since $-\frac{p}{2}+1-\delta$ and $-\frac{p}{2}+1+\eta+\delta$ belong to $(0,1)$ provided $0<\delta<\frac{1}{2m+1}$, we set them as $h$ for the first and the second term, respectively. 
Besides taking $a=\alpha$ and $l=p\nu<1$ which comes from (\ref{munu>0}), we can apply Lemma \ref{lm5.5} to the $\beta$-integrals and calculate (\ref{I2_u4}) as follows. 
\begin{align*}
&\hspace{-8pt}|I_2(Z^{-p},p)|\\
&=C\tau_-^{\frac{1}{2}-\eta}\int_{0}^{(\frac{t-r}{2}-k)_+}
\left(\frac{\alpha+2k}{k}\right)^{-1+\delta}(t-r-\alpha)^{m+1-(m+\frac{1}{2})p-\delta}\left(\log 3\frac{\alpha+2k}{k}\right)^{1-p\nu}
d\alpha\\
&\quad+C\tau_-^{m+\frac{5}{2}-(m+\frac{1}{2})p}\int_{\frac{t-r}{2}-k}^{t-r}
\left(\frac{\alpha+2k}{k}\right)^{-1-\eta-\delta}(t-r-\alpha)^{-1+\delta}\left(\log 3\frac{\alpha+2k}{k}\right)^{1-p\nu}d\alpha\\
&\le C\tau_-^{-\eta+m+\frac{3}{2}-(m+\frac{1}{2})p-\delta}(\log 3\tau_-)^{1-p\nu}\int_{0}^{(\frac{t-r}{2}-k)_+}
\left(\frac{\alpha+2k}{k}\right)^{-1+\delta}d\alpha\\
&\quad+C\tau_-^{-\eta+m+\frac{3}{2}-(m+\frac{1}{2})p-\delta}(\log 3\tau_-)^{1-p\nu}\int_{\frac{t-r}{2}-k}^{t-r}
(t-r-\alpha)^{-1+\delta}d\alpha\\
&\le C\tau_-^{-\eta}W(r,t)^{-1}E_1(T).
\end{align*}
\par
When $F(p,q;n)=0$ with $p=q$ or $F(p,q;n)>0$, we can treat them similarly, since $\phi_0$ and $Z(r,t)$ consist of $\tau_-$ only. 
Using Lemma \ref{lm5.4} for all of four, we can get the desired estimates. 
\noindent
\underline{The second line of (\ref{8})}
\par
Since $\phi_0$ and $W$ consist of $\tau_-$ only, we can calculate them by the same way as above and use (\ref{l}) for $I_2(W^{-q},q)$. Then, we have 
\begin{gather}
|I_2(\phi_0^{-q},q)|,\ 
|I_2(\phi_0^{1-q}W^{-1},q)|,\ 
|I_2(\phi_0^{-1}W^{1-q},q)|
\le C\tau_-^{-\eta+m+\frac{3}{2}-(m+\frac{1}{2})q},\\
|I_2(W^{-q},q)
\le C\tau_-^{-\eta+m+\frac{3}{2}-(m+\frac{1}{2})q}\times\left\{\begin{array}{ll}
\log \tau_-&\mbox{if}\quad F=0,\ p=q)\\
\tau_-^{1-\{(m+\frac{1}{2})p-(m+\frac{3}{2})\}q}&\mbox{Otherwise}.\end{array}\right.
\end{gather}
We shall use (\ref{5.13}) for the fist three estimates. 
As to the fourth, we can treat it by the same way for the corresponding case for $I_1$. 
Then we get the desired estimates and prove Proposition \ref{prop:12} for $I_2$. 
This completes the proof of Proposition \ref{prop:12}. 
\subsection{Estimates for $I_3$}
\par 
In this subsection, we assume $t-r>0$. 
We can calculated all the inequalities by the same way, we consider $I_3(\phi_0^{-p},p)$ for the case of $F=0$ with $p\neq q$ as an example. 
By (\ref{I3})and $\phi_0=\tau^{-m-\frac{1}{2}}$, we have 
\begin{align}
|I_3(\phi_0^{-p},p)|
&\le \int_{\frac{t-r-2k}{3}}^{t-r}
\left(\frac{t-r-\tau}{k}\right)^{m+\frac{1}{2}-(m-1)p-\eta}
\left(\frac{t-r-\tau+2k}{k}\right)^{-p}\nonumber\\
&\qquad\qquad\times
\left(\frac{\tau+\frac{\lambda_-}{2}+2k}{k}\right)^{-\frac{p}{2}}
\left(\frac{\tau-\frac{\lambda_-}{2}+2k}{k}\right)^{-(m+\frac{1}{2})p}d\tau.\label{I3_u1}
\end{align}
\par
By $t-r>0$, we shall divide it as follows. 
\begin{equation}
\label{divide3}
\int_{\frac{t-r-2k}{3}}^{t-r}d\tau
=\int_{\frac{t-r-2k}{3}}^{t-r-\frac{2}{3}k}d\tau
+\int_{t-r-\frac{2}{3}k}^{t-r}d\tau
=:A+B
\end{equation}
In the domain of integration of $A$, $t-r-\tau$ is equivalent to $t-r-\tau+2k$. 
By 
\[
\left(\frac{\tau+\frac{\lambda_-}{2}+2k}{k}\right)^{-\frac{p}{2}}\le C\tau_-^{-\eta}\left(\frac{t-r-\tau+2k}{k}\right)^{-\frac{p}{2}+\eta}, 
\]
we have 
\begin{align*}
A
&\le C\tau_-^{-\eta}\int_{\frac{t-r-2k}{3}}^{t-r-\frac{2}{3}k}
\left(\frac{t-r-\tau+2k}{k}\right)^{m+\frac{1}{2}-(m+\frac{1}{2})p}
\left(\frac{3\tau-(t-r)+4k}{k}\right)^{-(m+\frac{1}{2})p}d\tau.
\end{align*}
Moreover, dividing it as
\[
\int_{\frac{t-r-2k}{3}}^{t-r-\frac{2}{3}k}d\tau
=\int_{\frac{t-r-2k}{3}}^{\frac{t-r}{2}-\frac{2}{3}k}d\tau
+\int_{\frac{t-r}{2}-\frac{2}{3}k}^{t-r-\frac{2}{3}k}d\tau
=:A_1+A_2, 
\]
each integral can be calculated as follows. 
\begin{align*}
A_1
&\le C\tau_-^{-\eta+m+\frac{1}{2}-(m+\frac{1}{2})p}\int_{\frac{t-r-2k}{3}}^{\frac{t-r}{2}-\frac{2}{3}k}
\left(\frac{3\tau-(t-r)+4k}{2k}\right)^{-(m+1/2)p}d\tau\\
&\le C\tau_-^{-\eta+m+\frac{1}{2}-(m+\frac{1}{2})p}\\
&=C\tau_-^{-\eta}W(r,t)^{-1}.
\end{align*}
\begin{align*}
A_2
&\le C\tau_-^{-\eta-(m+\frac{1}{2})p}
\int_{\frac{t-r}{2}-\frac{2}{3}k}^{t-r-\frac{2}{3}k}\left(\frac{t-r-\tau+2k}{k}\right)^{m+\frac{1}{2}-(m+\frac{1}{2})p}d\tau\\
&\le C\tau_-^{-\eta-(m+\frac{1}{2})p}\\
&= C\tau_-^{-\eta}W(r,t)^{-1}.
\end{align*}
As to $B$ in (\ref{divide3}), we have  
\begin{align*}
B
&=
\int_{t-r-\frac{2}{3}k}^{t-r}
\left(\frac{t-r-\tau}{k}\right)^{m+\frac{1}{2}-(m-1)p-\eta}
\left(\frac{t-r-\tau+2k}{k}\right)^{-p}\\
&\qquad\qquad\qquad\times
\left(\frac{\tau+\frac{\lambda_-}{2}+2k}{k}\right)^{-\frac{p}{2}}
\left(\frac{\tau-\frac{\lambda_-}{2}+2k}{k}\right)^{-(m+\frac{1}{2})p}d\tau\\
&\le C\tau_-^{-\frac{p}{2}-(m+\frac{1}{2})p}\int_{t-r-\frac{2}{3}k}^{t-r}
\left(\frac{t-r-\tau}{k}\right)^{m+\frac{1}{2}-(m-1)p-\eta}d\tau. 
\end{align*}
Since 
\begin{equation}
\label{power}
m+\frac{3}{2}-(m-1)p-\eta>\frac{1}{2}+\frac{6}{2m+1}-\eta>0\qquad\mbox{for}\quad p<\frac{2m+5}{2m+1},\quad \eta=0\ \mbox{or}\ \frac{1}{2}, 
\end{equation}
we can integrate it as follows. 
\[
B
\le C\tau_-^{-\frac{p}{2}-(m+\frac{1}{2})p}
= C\tau_-^{-\eta}W(r,t)^{-1}.
\]
Therefore, we get 
\[
|I_3(\phi_0^{-p},p)|\le C\tau_-^{-\eta}W(r,t)^{-1}.
\]
\par
By the similar calculation, we have 
\begin{align}
|I_3(\phi_0^{1-p}Z^{-1},p)|
&\le C\tau_-^{-\eta+m+\frac{1}{2}-(m+\frac{1}{2})p}+C\tau_-^{-\eta-(m+\frac{1}{2})(p-1)-\frac{1}{p}}+C\tau_-^{-\frac{p}{2}-(m+\frac{1}{2})(p-1)-\frac{1}{p}},\nonumber\\
|I_3(Z^{-p},p)|
&\le C\tau_-^{-\eta+m+\frac{1}{2}-(m+\frac{1}{2})p}\log 3\frac{t-r+4k}{4k}+C\tau_-^{-\eta-1}+C\tau_-^{-\frac{p}{2}-1}.\label{I31_u3}
\end{align}
for the case of $F(p,q;n)=0$ with $p\neq q$, 
\begin{align*}
|I_3(\phi_0^{-p},p)|
&\le C\tau_-^{-\eta+m+\frac{1}{2}-(m+\frac{1}{2})p}+C\tau_-^{-\eta-(m+\frac{1}{2})p}+C\tau_-^{-\frac{p}{2}-(m+\frac{1}{2})p},\\
|I_3(\phi_0^{1-p}Z^{-1},p)|
&\le C\tau_-^{-\eta+m+\frac{1}{2}-(m+\frac{1}{2})p}+C\tau_-^{-\eta-(m+\frac{1}{2})(p-1)-\mu}+C\tau_-^{-\frac{p}{2}-(m+\frac{1}{2})(p-1)-\mu},\\
|I_3(Z^{-p},p)|
&\le C\tau_-^{-\eta+m+\frac{1}{2}-(m+\frac{1}{2})p}\times\left\{\begin{array}{cl}
\log\tau_-&\mbox{if}\quad F=0,\ p=q\\
\tau_-^{1-p\mu}&\mbox{otherwise}
\end{array}\right.\\
&\qquad\qquad+C\tau_-^{-\eta-p\mu}+C\tau_-^{-\frac{p}{2}-p\mu}.\label{I32_u3}
\end{align*}
for the case of $F(p,q;n)=0$ with $p=q$ or $F(p,q;n)>0$. 
Since it holds that 
\[
\log X\le \frac{1}{\delta}X^\delta\qquad(X\ge 1), 
\]
for any $\delta>0$, we shall apply it to (\ref{I31_u3}) with sufficiently small $\delta$. 
Thus, we have proved all the inequalities in the left column in Proposition \ref{prop:34} for all cases. 
\par
As for the right column in Proposition \ref{prop:34}, 
calculating them by the same way as above, we have 
\begin{align*}
|I_3(\phi_0^{-q},q)|
&\le C\tau_-^{-\eta+m+\frac{1}{2}-(m+\frac{1}{2})q}+C\tau_-^{-\eta-(m+\frac{1}{2})q}+C\tau_-^{-\frac{q}{2}-(m+\frac{1}{2})q},\\
|I_3(\phi_0^{1-q}W^{-1},q)|
&\le C\tau_-^{-\eta+m+\frac{1}{2}-(m+\frac{1}{2})q}+C\tau_-^{-\eta-(m+\frac{1}{2})(q-1)-\{(m+\frac{1}{2})p-(m+\frac{3}{2})\}}\\
&\qquad+C\tau_-^{-\frac{p}{2}-(m+\frac{1}{2})(q-1)-\{(m+\frac{1}{2})p-(m+\frac{3}{2})\}},\\
|I_3(W^{-q},q)|
&\le C\tau_-^{-\eta+m+\frac{1}{2}-(m+\frac{1}{2})q}\times\left\{\begin{array}{cl}
\log\tau_-&\mbox{if}\quad F=0,\ p=q\\
\tau_-^{1-\{(m+\frac{1}{2})p-(m+\frac{3}{2})\}q}&\mbox{otherwise}
\end{array}\right.\\
&\qquad\qquad+C\tau_-^{-\eta-\{(m+\frac{1}{2})p-(m+\frac{3}{2})\}q}+C\tau_-^{-\frac{p}{2}-\{(m+\frac{1}{2})p-(m+\frac{3}{2})\}q}.\label{I33_v3}
\end{align*}
Here we used (\ref{l}) for the third. 
By $m+\frac{5}{2}-(m+\frac{1}{2})p>0$ and (\ref{5.14}), the right column in Proposition \ref{prop:34} has been proved for all the cases. 
\subsection{Estimates for $I_4$}
\par
If $\Phi=\Phi(\tau_-)$, we have $|I_4(\Phi,p)|\le C\Phi(\tau_-)\w{I}_4$ by (\ref{I4}), where
\[
\w{I}_4:=
\int_{t-r}^t
\left(\frac{\tau+r-t}{k}\right)^{m+\frac{1}{2}-(m-1)p-\eta}
\left(\frac{\tau+r-t+2k}{k}\right)^{-p}
\left(\frac{2\tau+r-t+2k}{k}\right)^{-\frac{p}{2}}d\tau
\]
Since $W$ $Z$ and $\phi_0$ cosists of $\tau_-$ only, we shall consider $\w{I}_4$. 
\par
When $t-r+2k\ge t$, extending the upper limit, we have by (\ref{power}) 
\begin{align*}
\w{I}_4&\le C\tau_-^{-\frac{p}{2}}\int_{t-r}^{t-r+2k}\left(\frac{\tau+r-t}{k}\right)^{m+\frac{1}{2}-(m-1)p-\eta}d\tau
\le C\tau_-^{-\frac{p}{2}}.
\end{align*}
When $t-r+2k\le t$,we  divide the integral as 
\[
\int_{t-r}^td\tau=\int_{t-r}^{t-r+2k}d\tau+\int_{t-r+2k}^td\tau. 
\]
The first term is same as the one for $t-r+2k\le t$. 
By 
\[
\left(\frac{2\tau+r-t+2k}{k}\right)^{-\frac{p}{2}}\le \tau_-^{-\eta}\left(\frac{\tau+r-t}{k}\right)^{-\frac{p}{2}+\eta}
\]
the second term is calculated as below. 
\begin{align*}
\mbox{(the secont term)}
\le\tau_-^{-\eta}\int_{t-r+2k}^t\left(\frac{\tau+r-t}{k}\right)^{m+\frac{1}{2}-(m+\frac{1}{2})p}d\tau\le C\tau_-^{-\eta}.
\end{align*}
Therefore, we get 
\begin{align*}
|I_4(\phi_0^{-p},p)|
&\le C\tau_-^{-(m+\frac{1}{2})p-\frac{p}{2}}+C\tau_-^{-(m+\frac{1}{2})p-\eta}\\
|I_4(\phi_0^{1-p}Z^{-1},p)|
&\le C\tau_-^{-(m+\frac{1}{2})(p-1)-\frac{1}{p}-\frac{p}{2}}(\log 3\tau_-)^{-\nu}+
C\tau_-^{-(m+\frac{1}{2})(p-1)-\frac{1}{p}-\eta}(\log 3\tau_-)^{-\nu}\\
|I_4(\phi_0^{-1}Z^{1-p},p)|
&\le C\tau_-^{-1-\frac{p}{2}}(\log 3\tau_-)^{-p\nu}
+C\tau_-^{-1-\eta}(\log 3\tau_-)^{-p\nu}, 
\end{align*}
if $F(p,q;n)$ with $p\neq q$, 
\begin{align*}
|I_4(\phi_0^{-p},p)|
&\le C\tau_-^{-(m+\frac{1}{2})p-\frac{p}{2}}+C\tau_-^{-(m+\frac{1}{2})p-\eta}\\
|I_4(\phi_0^{1-p}Z^{-1},p)|
&\le C\tau_-^{-(m+\frac{1}{2})(p-1)-\mu-\frac{p}{2}}+
C\tau_-^{-(m+\frac{1}{2})(p-1)-\mu-\eta}\\
|I_4(\phi_0^{-1}Z^{1-p},p)|
&\le C\tau_-^{-p\mu-\frac{p}{2}}
+C\tau_-^{-p\mu-\eta}, 
\end{align*}
if $F(p,q;n)=0$ with $p=1$ or $F(p,q;n)>0$. 
Since $m+\frac{5}{2}-(m+\frac{1}{2})p>0$, all of them can be dominated by $C\tau_-^{-\eta}W(r,t)^{-1}$. 
Similary, we have 
\begin{align*}
|I_4(\phi_0^{-q},q)|
&\le C\tau_-^{-(m+\frac{1}{2})q-\frac{q}{2}}+C\tau_-^{-(m+\frac{1}{2})q-\eta}\\
|I_4(\phi_0^{1-q}W^{-1},q)|
&\le C\tau_-^{-(m+\frac{1}{2})(q-1)-\{(m+\frac{1}{2})p-(m+\frac{3}{2})\}-\frac{q}{2}}+
C\tau_-^{-(m+\frac{1}{2})(q-1)-\{(m+\frac{1}{2})p-(m+\frac{3}{2})\}-\eta}\\
|I_4(\phi_0^{-1}W^{1-q},q)|
&\le C\tau_-^{-\{(m+\frac{1}{2})p-(m+\frac{3}{2})\}q-\frac{q}{2}}
+C\tau_-^{-\{(m+\frac{1}{2})p-(m+\frac{3}{2})\}q-\eta}.
\end{align*}
By $m+\frac{5}{2}-(m+\frac{1}{2})q>0$ and (\ref{5.13}), the first two are dominated by $C\tau_-^{-\eta}W(r,t)^{-1}$. 
Moreover, using (\ref{5.14}), we have 
\begin{align*}
|I_4(\phi_0^{-1}W^{1-q},q)|
&\le C\tau_-^{-\eta-\mu+q(p-1)F}\\
&\le C\tau_-^{-\eta}\times\left\{\begin{array}{cl}
\tau_-^{\frac{1}{p}}(\log 3\tau_-)^{-\nu}(\log 3\tau_-)^{\nu}&\mbox{if}\quad F=0,\ p\neq q,\\
\tau_-^{-\mu}&\mbox{if}\quad F=0,\ p=q,\\
\tau_-^{-\mu+q(p-1)F}&\mbox{if}\quad F>0
\end{array}\right.\\
&\le C\tau_-^{-\eta}Z(r,t)^{-1}E_2(T)
\end{align*}
This completes the proof of Proposition \ref{prop:34}. 

\begin{center}
{\bf {\it Acknowledgements}}\\
\end{center}
\par
The author thanks Professor Hiroyuki Takamura (Tohoku Univ., Japan) for proposing this problem, 
stimulating discussions on the blow-up and local-existence of $L^p$ 
solutions and his valuable guidance. 


\end{document}